# A novel heuristic algorithm: adaptive and various learning-based algorithm

Sheng-Xue He

(Business School, University of Shanghai for Science and Technology, Shanghai 200093, China)

**Abstract**: A novel population-based heuristic algorithm called the adaptive and various learning-based algorithm (AVLA) is proposed for solving general optimization problems in this paper. The main idea of AVLA is inspired by the learning behaviors of individuals in a group, e.g. a school class. The algorithm formulates the following learning behaviors: a. Elite members will learn from each other; b. A common member will learn from some elite member and other common members; c. Members with unsatisfied performance will reflect their behavior after performance estimation; d. The whole group will reflect their behavior and try to improve if the performance of the group as a whole has not been improved for a long time. AVLA adopts the success-history based parameter adaptation to lighten the burden of parameter adjustment. To verify the efficiency of the AVLA, we apply it and its no-adaptation version with other eight well-known heuristics to 100 benchmark problems. The comparison clearly shows that AVLA performs as well as SHADE and the non-adaption version of AVLA outperforms others except AVLA and SHADE.

## 1. Introduction

In practice, heuristic algorithms have been used as a power tool to tackle all kinds of optimization problems for several decades. The reason behind the extensive use of heuristic algorithms is simple. On the one hand, a heuristic algorithm is usually easy to be applied due to its concise and clear implementation process and no-need of the information of derivatives. On the other hand, the real-life optimization problems generally have the non-smooth and nonlinear feature and are non-convex and constrained. Hundreds of heuristics have been proposed in past and many are very successful in practice. No Free Lunch Theory [1] explains the reason why so many heuristics are required. Researchers around the world have continuously made efforts to refine the existing algorithms and tried to create more effective ones in the past serval decades. After so many years of practice and accumulation, through it is still far from completely understanding the underlying operation mechanism of these algorithms, researchers have obtained many knowhows which can help them design more effective algorithms. In this paper, we will explore further on this way to design a novel effective algorithm based on the known rules and effective techniques. The new algorithm called the Adaptive and Various Learning-based Algorithm (AVLA) is inspired by the learning behaviors of individuals in a group. AVLA formulates the specified learning behaviors of individuals in divided subgroups at different learning stages. To improve the performance, the success history-based parameter adaptation is adopted in AVLA. Numerical analysis based on applying the new algorithm and other eight algorithms to 100 well-known benchmark test optimization problems demonstrates that the new algorithm can be a competitive and promising optimization tool for researcher and practitioners.

In the heuristic optimization field, researchers have proposed a huge number of algorithms in the past. Giving every algorithm a brief introduction is impossible in a paper due to the limited space.

Each one of the existing heuristic algorithms usually provides a feasible way of solving all kinds of optimization problems and performs well on some special type of problems. Heuristic algorithms can be grouped into four categories roughly. In the following, we will introduce only part of them especially the original metaheuristic algorithms not the various later-developed variants.

The first category includes all kinds of evolution algorithms. Many advanced and powerful and widely used heuristic algorithms come from this category. Genetic Algorithm (GA) [2] may be the first widely accepted evolution algorithm. Later, algorithms such as Differential Evolution (DE) [3] and the Covariance Matrix Adaptation Evolution Strategy(CMA-ES) [4] were introduced and successfully applied in practice. Tons of refined algorithms based on the earlier versions such as the CMA-ES including restarts with increasing population size [5] and the Combat Genetic Algorithm (CGA) [6, 7] were developed. Heuristic algorithms developed from the DE obtained huge success in practice due to their simple implementation processes and the capacity of obtaining satisfying optimization results. Algorithms such as JADE [8], Biogeography-Based Optimization (BBO) [9], the Success-History based parameter Adaptation for Differential Evolution (SHADE) [10], L-SHADE [11], the Backtracking Search Optimization Algorithm (BSO) [12] and Symbiotic Organism Search (SOS) [13] are just a few of them well-known to researchers.

The second category includes all kinds of swarm intelligent algorithms. One widely used swarm heuristic algorithm is the Ant Colony Optimization (ACO) presented by Dorigo in 1991 [14]. ACO is especially easy to be applied to discrete optimization problems. Dolphin Echolocation Optimization (DEO) proposed by Kaveh and Farhoudi in 2013 [15] is also fit for discrete optimization. Other two earlier developed and widely used swarm intelligent algorithms are Particle Swarm Optimization (PSO) [16] and Artificial Bee Colony algorithm (ABC) [17]. Animals and plants provide rich sources for researcher to design metaheuristic algorithms. Many of the algorithms achieve surprising success in a wide area. The following algorithms are inspired by the behaviors of all kinds of animals: Shuffled Frog-Leaping Algorithm (SFLA) [18], Honey-Bee Mating Optimization (HBMO) [19], Cuckoo Search with Levy Flights [20], Firefly Algorithm [21], Bat Algorithm [22] [23], Grey Wolf Optimizer (GWO) [24], Coral Reefs Optimization Algorithm (CRO) [25], Moth-flame optimization algorithm (MFO) [26], Ant Lion Optimizer (ALO) [27], Monkey King Evolution (MKE) [28], Whale Optimization Algorithm (WOA) [29], Dragonfly Algorithm [30], Grasshopper Optimization Algorithm (GOA) [31], Salp Swarm Algorithm (SSA) [32], Snake optimizer (SO) [33] and Giant Trevally Optimizer (GTO) [34]. In contrast to the algorithms inspired by animals' behaviors, algorithms such as Invasive Weed Optimization (IWO) [35], Plant Propagation Algorithm (PPA) [36] and Flower Pollination Algorithm (FPA) [37] are inspired by plant behaviors.

In addition to searching inspiring ideas from animals and plants, the physic and chemical phenomena and laws even math operations become the other rich area where many algorithms are inspired. Simulated annealing [38] may be the first very successful algorithm come from this category. The following are some representatives including Big Bung-Big Crunch algorithm (BBBC) [6], Central Force Optimization (CFO) [39], Gravity Search Algorithm (GSA) [40], Charged System Search (CSS) [41], Ray Optimization [42], Water Cycle Algorithm (WCA) [43], Black Hole Algorithm (BHA) [44], Global Mine Blast Algorithm (GMBA) [45, 46] [46], Radial Movement Optimization (RMO) [47], Enhanced Colliding Bodies Optimization (ECBO) [48], Stochastic Fractal Search (SFS) [49], Sine Cosine Algorithm (SCA) [50], Water Evaporation Optimization (WEO) [51], Multi-Verse Optimizer (MVO) [52], Vibrating Particles System algorithm (VPS) [53],

Tug of War Optimization (TWO) [54], Thermal Exchange Optimization (TEO) [55], Opposition-based tug of war optimization 2020 [56], and Equilibrium Optimizer (EO) [57].

The last category includes metaheuristic algorithms inspired by phenomena or actions associated with human society. Generally, Tabu Search [58] is viewed as the first one in this category. Harmony Search (HS) [59] is another algorithm in this category widely used in practice. Other algorithms such as Teaching-Learning-Based Optimization (TLBO) [60], Gaining-Sharing Knowledge based algorithm (GSK) [61], Life Choice-Based Optimizer (LCBO) (2019) [62] and Medalist Learning Algorithm (MLA) [63] [64] base their main idea on the learning activities of human being. Other representatives include Imperialist Competitive Algorithm (ICA) [65], League Championship Algorithm (LCA) [66], Jaya [67] and Levy flight Jaya algorithm (LJA) [68]. Among them, the implementation process of Jaya is very simple. A recently proposed algorithm in this category is Hiking Optimization Algorithm (HOA) [69]. The new algorithm-AVLA to be presented in this paper is also included in this category. AVLA is different from the other existing learning-based algorithm in this category by the concrete forms of formulating the learning behaviors.

The main contribution of this paper is presenting a novel powerful metaheuristic algorithm called the Adaptive and Various Learning Algorithm (AVLA). By comparing it with other well-known or newly proposed heuristic algorithms, we make it clear that AVLA is a competitive optimization tool for researchers and practitioners.

The rest of the paper is organized as follows. Section 2 introduces the main operations and the application process of AVLA. The optimization mechanism is also explained in Section 2 by focusing on the realization of balancing exploitation and exploration during the optimization. Section 3 verifies the efficiency of AVLA by implementing it and nine other algorithms to solving a series of optimization problems including 29 well-used test problems, CEC2014 test suite, CEC 15 learning-based benchmark suite, CEC2022 test suits and 14 real engineering optimization problems. Section 4 summarizes the paper and discusses several promising research directions.

## 2. Adaptive and Various Learning-based Algorithm

### 2.1 Basic ideas of AVLA

The individual learning behaviors in a group are complex. In general, members will choose their learning strategies according to their relative positions or performances in the group. A member with high performance is more likely to learning from other members with high performance. A member with common or unsatisfied performance is prone to emulate the learning strategy of some other individual with high performance. And at the same time common learners are easily influenced by their pals or near friends. After the performances are estimated and all members are sorted according to their performances, the members with unsatisfied performance may conduct some form of reflections to refine their scores. Sometimes if the performance cannot be refined, the member with unsatisfied performance will be replaced by a new member. If the performance of the group as a whole has not been improved for a long time, all the members may need to conduct some type of reflections about their learning strategies and expect to refine the performance of the whole group. The behaviors described above are natural and reasonable for most of the learning groups composed of human being.

In this paper, we will model the above learning actions into a novel heuristic algorithm called the adaptive and various learning-based algorithm (AVLA). AVLA is a population-based meta-heuristic algorithm for global and constrained optimization. A member in a learning group stands for a solution to the optimization problem in question. The learning group is the population of solutions. The performance of a member is related to the objective function value. When constrained optimization is considered, the performance is corresponding to the objective function value with penalty on the violation of constraints.

**2.2 Designation of AVLA**

## 2.2.1 Problem description

At first, let us introduce the main notions to be used to describe our optimization problem. Let $X \coloneqq \{x_i | i = 1,2, \ldots, N\}$ be the solution set with $x_i$ as its representative. As mentioned above, a solution is a member in the learning group. Here $N$ is the size of the population (or called group in this paper). A solution is composed of $D$ elements. Let $x_i = (x_{i,1}, x_{i,2}, \ldots, x_{i,D})$. We assume that any element needs to be chosen from a given set. That is $x_{i,j} \in V_j$ for all $i$ where $V_j$ is the feasible range or set of the $j$th element of any member. Assume that $\alpha_j \leq x_{i,j} \leq \beta_j$, $\forall i$. The objective function value (or the penalty function value corresponding to a constrained optimization problem) with respect to member $x_i$ is denoted by $fit(x_i)$. We aim to minimizing the objective function value. The initial solutions are generated randomly in the feasible space. For a member $x_i$, its initial $j$th element $x_{i,j}$ is generated as follows:

$$x_{i,j} = \alpha_j + rand[0,1](\beta_j - \alpha_j) \tag{1}$$

where $rand[0,1]$ is random number picked from interval $[0,1]$ according to the uniform distribution.

## 2.2.2 Learning process modeling

In this subsection we will model the learning actions of members in a learning group. Firstly, we will deal with the learning behaviors of members with high performance. We will call a member with high performance an elite in this paper. Secondly, the learning actions of common learners will be formulated into mathematic functions. At the end, we will model the reflection actions of members with unsatisfied performance and the whole group, respectively.

*a. Learning among elites*

Confucius said: "Let me take a stroll with any two people, and I can always be sure of learning something from them. I can take their good points and emulate them, and I can take their bad points and correct them in myself." (In Chinese, 子曰: 三人行, 必有我師焉; 擇其善者而從之, 其不善者而改之.) The elites in a learning group can put the above saying in practice as follows. At first, an elite can choose two other elites randomly as his/her models to emulate their strategies; then by comparing performances, the elite can determine whether to follow or to keep away from the selected elites.

In view of that the actual learning efficiency is always influenced by all kinds of factors and shows obvious uncertainty, we assume that the learning of an individual will be realized through three steps. Three steps are the ideal learning, practical learning, and actual learning. They are modeled in equations (2), (3) and (4), respectively, as given below:

$$v_{e,j}^t = x_{e,j}^t + s_1 F_e \left( x_{e1,j}^t - x_{e,j}^t \right) + s_2 F_e \left( x_{e2,j}^t - x_{e,j}^t \right) \tag{2}$$

$$v_{e,j}^{t+1} = \begin{cases} v_{e,j}^t & \text{if } rand[0,1] \leq CR_e \text{ or } j = j_{rand} \\ x_{e,j}^t & \text{otherwise} \end{cases} \quad (3)$$

$$x_e^{t+1} = \begin{cases} v_e^{t+1} & \text{if } fit(v_e^{t+1}) < fit(x_e^t) \\ x_e^t & \text{otherwise} \end{cases} \quad (4)$$

In above, $e$ indicates the elite in question; $e1$ and $e2$ are two randomly chosen different elites. $t$ is the learning stage that is corresponding to the iteration of AVLA. $s_1$ and $s_2$ are two sign parameters. If $fit(x_e^t) > fit(x_{e1}^t)$, $s_1$ equals 1; or else, $s_1$ equals -1. Similarly, if $fit(x_e^t) > fit(x_{e2}^t)$, $s_2$ equals 1; or else, $s_2$ equals -1. $F_e$ is *the learning acceptation rate* of member $e$. $F_e$ will be limited to [0, 1]. The rule of how to change $F_e$ will be introduced in **Section 2.2.3**. $v_{e,j}^t$ can be viewed as the ideal result from learning from other two elites because in this situation, $e$ can refine all his/her elements in a full scale.

In reality, it is nearly impossible for a member to learn everything from another member. A practical learning is to just emulate part of the features of the selected model. Here a feature stands for an element of a solution to the optimization problem in question. Equation (3) formulates the above thought with a stochastic judgement. $j_{rand}$ is an index chosen randomly from 1, 2, …, D. $rand[0,1]$ is a random number chosen from interval $[0,1]$ according to uniform distribution. $CR_e$ is *the practical learning factor*. The rule of how to change $CR_e$ will be introduced in **Section 2.2.3**.

We assume that an elite only accepts the change which can improve his/her current performance. Equation (4) models this decision-making behavior called the actual learning. The result from the above three steps will be the new state of elite $e$ to be estimated later.

*b. Learning by commons*

Confucius said: "I listen a lot, pick the best of it, and follow it; I observe a lot and take note of it. This is the best way for me to learn." (In Chinese, 子曰: 多聞, 擇其善者而從之; 多見而識之, 知之次也.) For a common learner whose performance is lower than any elite, the learning process is similar to the elites but with different learning models to emulate. A common member stands for an individual other than elites in AVLA. We assume that a common member may choose to learn from an elite and one common member or decide to learn from two common members. The learning process followed by a common member consists of three steps similar to that related to the learning process among elites. The three steps composed of the ideal learning, practical learning and actual learning are formulated in equations (5), (6) and (7), respectively, as follows:

$$v_{i,j}^t = \begin{cases} x_{i,j}^t + s_1 F_i(x_{i1,j}^t - x_{i,j}^t) + s_2 F_i(x_{i2,j}^t - x_{i,j}^t) & \text{if } rand[0,1] > \text{LE}(t) \\ x_{i,j}^t + F_i(x_{e,j}^t - x_{i,j}^t) + s_2 F_i(x_{i2,j}^t - x_{i,j}^t) & \text{otherwise} \end{cases} \quad (5)$$

$$v_{i,j}^{t+1} = \begin{cases} v_{i,j}^t & \text{if } rand[0,1) \leq CR_i \text{ or } j = j_{rand} \\ x_{i,j}^t & \text{otherwise} \end{cases} \quad (6)$$

$$x_i^{t+1} = \begin{cases} v_i^{t+1} & \text{if } fit(v_i^{t+1}) < fit(x_i^t) \\ x_i^t & \text{otherwise} \end{cases} \quad (7)$$

In equation (5), $i$ denotes the common learner in question; $e$ stands for a randomly selected elite; $i1$ and $i2$ are two members randomly chosen from the common members different from $i$ and $e$. $s_1$ and $s_2$ are two sign parameters. If $fit(x_i^t) > fit(x_{i1}^t)$, $s_1$ equals 1; or else, $s_1$ equals -1. Similarly, if $fit(x_i^t) > fit(x_{i2}^t)$, $s_2$ equals 1; or else, $s_2$ equals -1. $F_i$ and $CR_i$ are the learning acceptation rate and the practical learning factor of common member $i$. The explanation of equations (6)-(7) is similar to the explanation of equations (3)-(4).

In equation (5), LE(t) stands for the probability of that a common member decides to learn from an elite not another common member. LE(t) is defined by the following equation:

$$LE(t) = \left[1 + e^{\frac{2\gamma}{maxNumIter}\left(\frac{maxNumIter}{2}-t\right)}\right]^{-1} \quad (8)$$

In equation (8), $maxNumIter$ is the maximal number of iterations allowed to run AVLA once to solve an optimization problem. $t$ is the counter of the current iteration. $\gamma$ is called the truncated radius of learning curve. In this paper, we just assign 6 to $\gamma$. This curve function (8) is also used in MLA [64]. With the increase of $t$, the value of $LE(t)$ gradually increases from a low positive value near to zero to a high value a little bit less than 1. With LE(t) in (5), the common member is more likely to learn from other common learners at the beginning, but with the increase of iteration, he/she will be prone to learning from some elite.

*c. Reflection after sorting*

Confucius said: "Learning without thinking leads to perplexity. Thinking without learning leads to trouble." (In Chinese, 子曰：學而不思則罔，思而不學則殆.) After the previous mentioned learning, members will be sorted from the best to the worst according to their current performances. We assume that members will conduct some forms of reflections to further refine their performance after the sorting.

Before presenting the concrete reflection way, we need the following notation. Let $x_{i,j}$ that is the $j$th element of the member $i$ be limited to the interval $[\alpha_j, \beta_j]$. We define $x_{i,j}^R$ the reflection position of $x_{i,j}$ in $[\alpha_j, \beta_j]$ such that:

$$x_{i,j}^R = \alpha_j + \beta_j - x_{i,j} \quad (9)$$

All the reflection positions of the elements of $x_i$ form a new spot in the searching space denoted by $x_i^R$. We call $x_i^R$ the *opposite position* of $x_i$ in the search space.

Now with the new conception of opposite position, we can formulate the reflection behavior of a member with unsatisfied performance as follows:

$$x_i = \begin{cases} x_i^R & \text{if } fit(x_i^R) < fit(x_i) \\ x_i^{STOCH} & \text{otherwise} \end{cases} \quad \forall \, N - n_E < i \leq N \quad (10)$$

In equation (10), $N$ and $n_E$ are the size of population and the size of elites, respectively. $x_i^{STOCH}$ stands for a solution generated randomly with equation (1). Equation (10) means that if a member's position is at the tail part of the list of sorted members, he/she will do reflection of trying his/her opposite position in the search space. If his/her opposite position has a better outcome, he/she will occupy the new position; or else, he/she will choose a random position, in other words being replaced by a new member.

In addition to the above reflection carrying out by members with unsatisfied performance, the whole group may conduct some reflection if the performance of the group as a whole is not improved for a long time. In AVLA, the performance of the group as a whole is indicated by the best-so-far objective function value. If the best-so-far objective function value has not been improved for a given number $n_R$ of successive iterations, all the members will conduct reflection as follows:

$$x_i = \begin{cases} x_i^R & \text{if } fit(x_i^R) < fit(x_i) \text{ and } i \leq N - n_E \\ x_i & \text{if } fit(x_i^R) \geq fit(x_i) \text{ and } i \leq N - n_E \\ x_i^R & \text{if } N - n_E < i \leq N \end{cases} \quad (11)$$

Equation (11) means that if a member who is not a member with unsatisfied performance, he/she will check his/her opposite position and occupy it if it can lead to better performance than his/her current position. But if the member has an unsatisfied performance, he/she will be forced to move to his/her opposite position.

### 2.2.3 Parameter adaptation

To obtain the expected performance of a heuristic algorithm applied to a given problem, researchers generally need to adjust the parameters by trial-and-error. In general, there is no fixed parameter setting suitable for all problems and even the different phases during the process of optimization. To lighten the above burden and improve the resilience of a heuristic algorithm, parameter adaptation is adopted in many algorithms, e.g. JADE. The successful applications in the past in optimization field encourage the author to adopt it in AVLA. The adaptive method adopted in AVLA is the success-history based parameter adaptation first used in SHADE.

Table 1 Historical memory of the learning acceptance rate and the practical learning factor

| Index | 1 | 2 | ... | H-1 | H |
|---|---|---|---|---|---|
| $M_{CR}$ | $M_{CR,1}$ | $M_{CR,2}$ | ... | $M_{CR,H-1}$ | $M_{CR,H}$ |
| $M_F$ | $M_{F,1}$ | $M_{F,2}$ | ... | $M_{F,H-1}$ | $M_{F,H}$ |

AVLA maintains a historical memory with $H$ entries for both the learning acceptance rate $F$ and the practical learning factor $CR$ as shown in Table 1. The initial values of all the entries will be set to 0.5. At the beginning of each iteration of AVLA, every member $i$ in the learning group will select his/her learning acceptance rate $F_i$ and practical learning factor $CR_i$ as follows:

$$CR_i = randn_i(M_{CR,r_i}, 0.1) \qquad (12)$$

$$F_i = randc_i(M_{F,r_i}, 0.1) \qquad (13)$$

In equations (12) and (13), $r_i$ is a randomly selected index from $[1, H]$. $randn_i(M_{CR,r_i}, 0.1)$ means to generate $CR_i$ according to a normal distribution with mean $M_{CR,r_i}$ and standard deviation 0.1. If the generated $CR_i$ is outside of $[0,1]$, it will be truncated. $randc_i(M_{F,r_i}, 0.1)$ means to generate $F_i$ according to a Cauchy distribution with location parameter $M_{F,r_i}$ and scale parameter 0.1. If the generated $F_i$ is greater than 1, it will be replaced by 1; if the value is less than 0, the generation process using Cauchy distribution will be repeated until a proper value of $F_i$ is generated.

At the end of each iteration, the historical memory will be updated. During the learning process in an iteration, if a member $i$ successfully improves his/her performance, his/her learning acceptance rate $F_i$ and his/her practical learning factor $CR_i$ will be stored in sets $S_{CR}$ and $S_F$, respectively. The historical memory is updated as follows:

$$M_{CR,k,t+1} = \begin{cases} mean_{WL}(S_{CR}) & if\ S_{CR} \neq \emptyset \\ M_{CR,k,t} & otherwise \end{cases} \qquad (14)$$

$$M_{F,k,t+1} = \begin{cases} mean_{WL}(S_F) & if\ S_F \neq \emptyset \\ M_{F,k,t} & otherwise \end{cases} \qquad (15)$$

In equations (14) and (15), the subscript $k$ is the index of the memory entry to be updated. The initial value of $k$ is set to 1. Then whenever a new element is inserted into the memory, the value of $k$ will be increased by 1. If the value of $k$ becomes greater than $H$, it will be reset to 1.

The equations (14) and (15) use the weighted Lehmer means defined in equations (15) and (16):

$$mean_{WL}(S_{CR}) = \frac{\sum_{k=1}^{|S_{CR}|} w_k S_{CR,k}^2}{\sum_{k=1}^{|S_{CR}|} w_k S_{CR,k}} \qquad (16)$$

$$mean_{WL}(S_F) = \frac{\sum_{k=1}^{|S_F|} w_k S_{F,k}^2}{\sum_{k=1}^{|S_F|} w_k S_{F,k}} \qquad (17)$$

The weight $w_k$ is defined as follows:

$$w_k = \frac{\Delta fit_k}{\sum_{l=1}^{|S_{CR}|} fit_l} \tag{18}$$

$$\Delta fit_k = |fit(v_k^t) - fit(x_k^t)| \tag{19}$$

In equation (19), $\Delta fit_k$ is the absolute value of the change of the objective function values for member $k$ after the three-steps individual learning.

### 2.2.4 Increase the size of elites

The size of elites has a noticeable impact on the convergence of AVLA. At the earlier phase of search, the information held by elites is very likely uncertain and misleading regarding to the final best solution. Learning from them at the earlier phase of the search will do no obvious good to other members. In view of the above observation, only a few of members with relative high performance should be treated with as elites at the earlier phase of learning. When the optimization keeps going on, elites as a whole are more likely holding useful information of the direction of further refinement. At the same time, learning among elites can speed up the convergence to the promising position of the search space at the ending phase of learning. In view of the above observation, we will adapt the size of elites according to the following equation:

$$n_E(t) = round\left[3 + \frac{t(0.2N-3)}{maxNumIter}\right] \tag{20}$$

In equation (20), $t$ is the index of the current iteration and $maxNumIter$ is the maximal number of iterations predetermined for AVLA. It is easy to see that the initial size of elites is set to 3. With the increase of iterations, the size of elites will gradually be increased to $0.2N$ at the end.

The initial and final size of elites may be changed with respect to the problem in question. The above choice is just based on the experience of the author and it works most of time.

### 2.2.5 The process of the Algorithm

In the following we sum up the implementation process of AVLA as given previously.

**Step 1 Initialization**

Initialize values for $N$, $maxNumIter$, $n_R$ and $H$. And set $nR = 0$.

**Step 2 Generate the initial learning group**

Generate the initial learning group $X \coloneqq \{x_i | i = 1,2,...,N\}$ using equation (1). And sort the members according to their performances.

**Step 3 Conduct learning**

    **Step 3.1** For every elite member, carry out the three-steps learning by using equations (2-4).

    **Step 3.2** For every common learner, carry out the three-steps learning by using equations (5-7).

    **Step 3.3** Sort the members according to their new performances.

**Step 4 Update historical memory**

Update the historical memory by equations (14) and (15).

**Step 5 Do necessary reflections**

If the best-so-far solution has not changed after $n_R$ successive iterations, in other words $nR$ equals $n_R$, set $nR = 0$ and execute **Step 5.1**; or else, increase $nR$ by 1 and carry out **Step 5.2**.

    **Step 5.1** Do reflection of the whole group using equation (11). Jump to **Step 5.3**.

    **Step 5.2** Do reflection by members with unsatisfied performance using equation (10).

    **Step 5.3** Sort the members according to their new relative performance.

**Step 6 Check terminal condition**

If the terminal criteria are satisfied, output the best-so-far solution and end the program; or else, go to **Step 3**.

## 2.3 Analysis of the searching mechanism

An effective heuristic algorithm needs to properly balance its exploitation and exploration during the search. Exploitation means to search around some obtained solutions which usually have relatively high performance; exploration means to search more widely in the feasible search space. Generally speaking, meta-heuristic algorithm needs to stress on exploration in the first half of the search, but make more efforts on exploitation in the second half. The design of AVLA reflects the above idea.

By dividing the members into subgroups and then carrying out different learning strategies, we can gradually strengthen the exploitation by learning among elites with the gradually increasing size of elites and dynamically balance the exploitation and exploration by the learning of common learners with the increasing probability of learning from an elite not another common learner by a common member.

A typical individual learning behavior includes three steps. The ideal learning step determines the possible improvement quantity for every element of the solution. Realizing these improvement quantities will lead the learner nearer to his/her model in all dimensions. To avoid the premature convergence induced by the overlearning of the first step, the second learning step only accepts part of the possible improvement quantities. The last step of individual learning accepts the new solution with actual improvement. The last step guarantees the learning group as a whole being improved continuously.

When carrying out the ideal learning of a common member, we use a parameter to control the probability of choosing an elite not a common member to learn. The probability will gradually increase as the search keeps going on. This design is to avoid the ineffective learning of common members at the earlier stage and to put more stress on emulating promising elites at the ending stage.

It is obvious that the above learning emphasizes the probable improvement along the highly promising direction. The stochastic feature still needs to be strengthened in some way. To strengthen the stochastic feature, in other words the exploration capacity of AVLA, the reflection operations are required. Two types of reflections are used in AVLA to strengthen the exploration.

The reflection of a member with unsatisfied performance is realized by checking the opposite position, occupying the opposite position with better performance or a random position in the search space. This type of reflections tries to more evenly explore the search space and at the same time strengthens the random feature of AVLA.

Different from the above reflection conducted only by part of the members, the second type of reflections will be conducted by all members and be carried out only when the best-so-far solution has not been improved for a given number of successive iterations of AVLA. This can be viewed as a trial to restart the search with all members to be relocated to more promising positions. During the whole group reflection, if a member's performance is not unsatisfied, he/she will choose the better position between his/her current position and the corresponding opposite position. To strengthen the randomness of search, AVLA will force the members with unsatisfied performance to move to their opposite positions.

As mentioned in **Section 2.2.4**, the operation of increasing the size of elites is to balance the exploration and exploitation during the whole search. The gradually increased size of elites makes

AVLA to emphasize exploring more wide area at the beginning and then gradually put more energy to exploit the nearby area about the best-so-far solutions.

The adoption of the adaptive control of parameters with success history is based on the following understandings – "There is no fixed control parameter setting suitable for various problems or even at different optimization stages of a single problem." [8].

## 3. Numerical Experiments

### 3.1 Test on 100 problems

In this section, we will apply AVLA, the various-learning algorithm (VLA) and other eight heuristic algorithms to various optimization problems and try to figure out the relative performance of AVLA among them.

VLA is the variant of AVLA without the parameter adaptation. In VLA, the learning acceptation rate $F_i, i = 1,2,..,N$ is replaced by a random number $\text{Rand}_i[0,1]$ generated according to the uniform distribution on interval [0,1]. The practical learning rate CR is replaced by a given number which is in interval (0,1). In VLA, the operations associated with historical memory will be omitted.

Since the total number of the existing heuristic algorithms is large, we can only choose some representatives whose great performances have been proven and accepted by many in practice. The other eight algorithms to be compared include CMA-ES[7], SHANDE[10], TLBO[60], LCBO[62], GMBA[46], EO[57], MFO[26], and SO[33]. Among them, CMA-ES and SHADE are evolution algorithms. SHADE is the winner of CEC 2014 and its various variants keep the winning records in the following CEC series. TLBO and LCBO are two grouped into the human society-inspired algorithms. GMBA and EO come from the heuristics inspired by physical or chemical or mathematical laws. MFO and SO are two swarm intelligent heuristic algorithms. In addition to their even distribution across the categories, their birth times span decades. CMA-ES first appeared in 2005, but SO is recently presented in 2022.

The parameters used in these algorithms are summarized in Table 2. The size of population and the maximal number of iterations are set to 50 and 2000, respectively. Except CMA-ES, all others will run exact 2000 iterations before being terminated. Due to its inner mechanism, CMA-ES may terminate earlier. The other parameters are usually set as the authors suggested in their papers. When solving a problem, an algorithm will be run 30 times independently. We will calculate the mean and standard deviation (std) of the best objective function values and note down the best one among the 30 results. The author thinks the size of population and maximal number of iterations set as above is proper to show these algorithm relative performances. Note that a fixed setting of parameters for a given algorithm cannot be optimal for all problems. But supplying each problem with an optimal setting of parameters is just beyond the scope of this paper.

The codes of CMA-ES and SHADE comes from the corresponding authors' websites. Other algorithms are coded by the author in Java, but it is easy for reader to find avail codes in MATLAB or Python in some open source software, e.g. MEALPY [70]. In other words, the results given in this paper can be easily verified by interested readers. The following fact is worth addressing here. Though the results of a problem from different runs of the same algorithm are usually different due to the randomness of the algorithm, the statistic mean, the standard deviation and the best objective value seldom change or change a little from the 30-times runs. For example, if two results are 4.00 and 4.01 associated with two different algorithms, respectively, the difference 0.01 between them may look trivial, but this difference usually will not disappear or change its sign from positive to negative after many tries 0f 30-times runs.

Table 2 Algorithms to be used and their parameters.

| Algorithm | Parameters |
|---|---|
| CMA-ES | $N$=50, $\mu = 30$, $maxNumIter$=2000, others are set with default values given by the authors |
| SHADE | $N$=50, $maxNumIter$=2000, $p_{min} = 0.05$, $H$=50 |
| TLBO | $N$=50, $maxNumIter$=2000 |
| LCBO | $N$=50, $maxNumIter$=2000, other parameters are set as the authors suggested in their paper |
| GMBA | $N$=50, $maxNumIter$=2000, the number of explorations equals 1000, the reduction constant is 120 |
| EO | $N$=50, $maxNumIter$=2000, the size of memory is 4 |
| MFO | $N$=50, $maxNumIter$=2000, the parameter of logarithmic spiral is 1.0 |
| SO | $N$=50, $maxNumIter$=2000. |
| VLA | $N$=50, $maxNumIter$=2000, $CR$=0.25, $n_R$=10. |
| AVLA | $N$=50, $maxNumIter$=2000, $H$=50, $n_R$=6. |

Since our goal is to find out if AVLA is the best algorithm or the second best, we will use a simple way to compare the results. For a given group of problems, we will highlight in red color the resulted means which are at the first place for any problem in the group. Then we count the number of highlighted results which are at the first place for each algorithm. At last we compare the obtained number following the rule that the bigger the number of first places is, the better the corresponding algorithm performs on the given group of problems. The above method is easy and fit for our aim. Readers can check the correctness easily using the given data.

The first group of optimization problems includes 29 well-known test problems which are widely used in tons of papers. These problems are summarized in the **Appendix**. Different from the common use of these problems as a whole, we divide them into 4 sub-groups and check the performances of algorithms on the separated subgroups. In doing so, we can easily figure out some usually overlooked features of the algorithms in question.

The first sub-group includes 7 unimodal benchmark functions. The related results are shown in Table 3. In the following tables, the short hands "Ave", "Std" and "Best" are used to indicate the mean, the standard deviation and the best objective function value resulted from 30 times independent runs of an algorithm for solving a specified problem. "Num of Tops" is used to indicate the number of an algorithm's means occupying the first place. The numbers in the row with title "Rank" indicate the relative rank of an algorithm based on the number of the algorithm's means occupying the first place. The data in Table 3 show that LCBO and SO achieve the first and second in the rank. Both AVLA and SHADE only obtain 2 first places with the resulted means and their ranks are 4. VLA performs even worse than AVLA and only obtains 1 outstanding result. But if we look from a different view point, the better performance on this subgroup to some extent indicates the stronger local searching power of the algorithm. This stronger exploitation capacity of an algorithm may harm its performance on other types of problems. The above conjecture will be proven by the following experiments.

There are 6 multimodal benchmark functions included in the second subgroup. The results of the problems in this subgroup obtained by 10 algorithms are given in Table 4. AVLA is the best algorithm for this subgroup. SHADE and SO take the second and third places in the rank, respectively. The performance of VLA is mediocre this time just as it does for the first subgroup.

When applying 10 algorithms to the 10 multimodal problems with fixed dimension in the third subgroup, SHADE rises to the first place in the rank. VLA and AVLA tie at the second rank. Except the problem F15, the above three algorithms all obtain the final best means for all other problems. EO has an impressive performance this time by obtaining 7 best means. The data associated with the above analysis are presented in Table 5.

The fourth subgroup includes 6 composite benchmark functions. When applied to solving the problems in this subgroup, the performance of SHADE is disappointing. In contrast, the performances of both VLA and AVLA are outstanding. The first rank is commonly occupied by VLA and AVLA. MFO gains the second rank with 2 best means. VLA works out the best for problems F26 and F29 and AVLA obtains the best for F25 and F28. The above difference demonstrates the respective values of AVLA and VLA. The data associated with the above analysis are given in Table 5.

If we put together the rank information of 10 algorithms for the 29 problems in the first group, we obtain Table 7. The data in Table 7 show that AVLA, SHADE and VLA occupy the first, second and third places in the rank, respectively, with 19, 16 and 13 best final means resulted from 30 times runs of each algorithm on the 29 test problems.

**Table 3** Optimization results and comparison for unimodal benchmark functions.

| Fun. | Index | CMA-ES | SHADE | TLBO | LCBO | GMBA | EO | MFO | SO | VLA | AVLA |
|---|---|---|---|---|---|---|---|---|---|---|---|
| F1 | Ave | 1.13E04 | 1.79E-82 | **0** | **0** | 1.29E-03 | 1.65E-278 | 2.04E-69 | **0** | 4.4E-33 | 5.71E-84 |
|  | Std | 9E03 | 6.31E-82 | 0 | 0 | 5.92E-03 | 0 | 1.12E-68 | 0 | 8.68E-32 | 3.12E-83 |
|  | Best | 1.07E-17 | 4.84E-88 | 0 | 0 | 6.36E-13 | 4.95E-289 | 2.27E-109 | 0 | 6.02E-36 | 0 |
| F2 | Ave | 2.23E01 | 3.46E-45 | 3.66E-290 | **0** | 1.6E-01 | 8.41E-153 | 3.33E-01 | 1.04E-315 | 2.75E-26 | 1.55E-28 |
|  | Std | 1.38E01 | 9.58E-45 | 0E00 | 0 | 5.56E-01 | 2.1E-152 | 1.83E00 | 0E00 | 4.13E-26 | 1.57E-28 |
|  | Best | 7.36E-10 | 8.74E-49 | 1.07E-303 | 0 | 2.93E-07 | 9.86E-156 | 4.42E-54 | 4.9E-324 | 1.39E-27 | 0 |
| F3 | Ave | 4.71E05 | 1.22E-81 | **0** | **0** | 1.77E-03 | 2.39E-278 | 3.33E02 | **0** | 2.17E-30 | **0** |
|  | Std | 3.31E05 | 2.63E-81 | 0 | 0 | 4.98E-03 | 0E00 | 1.83E03 | 0 | 5.49E-30 | 0 |
|  | Best | 1.85E-17 | 1.21E-85 | 0 | 0 | 8.25E-09 | 1E-284 | 3.7E-109 | 0 | 1.3E-33 | 0 |
| F4 | Ave | 1E02 | 5.25E-35 | 1.53E-241 | **0** | 9.37E-04 | 2.61E-98 | 1.86E-09 | 4.13E-309 | 1.1E-01 | 9.16E-22 |
|  | Std | 0E00 | 1.1E-34 | 0E00 | 0 | 2.28E-03 | 1.41E-97 | 4.27E-09 | 0E00 | 9.46E-02 | 3.82E-21 |
|  | Best | 1E02 | 2.53E-37 | 1.51E-254 | 0 | 6.33E-07 | 7.34E-106 | 1.69E-13 | 2.36E-314 | 8.25E-03 | 0 |
| F5 | Ave | 8.97E07 | **0** | 6.69E00 | 3.9E00 | 3.28E01 | 4.97E00 | 9.15E03 | 8.87E00 | 1.81E00 | 8.65E-29 |
|  | Std | 7.16E07 | 0 | 1.02E00 | 3.44E-01 | 6.91E01 | 8.71E-01 | 2.74E04 | 9.04E-02 | 1.82E00 | 3.58E-28 |
|  | Best | 5.5E-19 | 0 | 3.3E00 | 3.12E00 | 2.52E-01 | 3.7E00 | 5.08E-02 | 8.64E00 | 2.23E-03 | 0 |
| F6 | Ave | 4.33E03 | **0** | **0** | **0** | 1.2E00 | **0** | 3E-01 | **0** | **0** | **0** |
|  | Std | 5.68E03 | 0 | 0 | 0 | 3.96E00 | 0 | 7.94E-01 | 0 | 0 | 0 |
|  | Best | 0 | 0 | 0 | 0 | 0 | 0 | 0E00 | 0 | 0 | 0 |
| F7 | Ave | 1.01E01 | 5.4E-04 | 1.84E-04 | 3.29E-05 | 1.14E-02 | 5.97E-05 | 4.49E-03 | **2.05E-05** | 3.45E-05 | 3.46E-04 |
|  | Std | 1.04E01 | 1.99E-04 | 1.38E-04 | 1.99E-05 | 3.03E-02 | 4.68E-05 | 4.1E-03 | 1.72E-05 | 2.35E-05 | 1.53E-04 |
|  | Best | 7.7S-03 | 2.3E-04 | 2.69E-05 | 2.7E-06 | 2.73E-04 | 4.73E-06 | 5.29E-04 | 4.15E-07 | 2.48E-06 | 9.72E-05 |
| Num Of Tops |  | 0 | 2 | 3 | **5** | 0 | 1 | 0 | 4 | 1 | 2 |
| Rank |  | 6 | 4 | 3 | **1** | 6 | 5 | 6 | 2 | 5 | 4 |

**Table 4** Optimization results and comparison for multimodal benchmark functions.

| Fun. | Index | CMA-ES | SHADE | TLBO | LCBO | GMBA | EO | MFO | SO | VLA | AVLA |
|---|---|---|---|---|---|---|---|---|---|---|---|
| F8 | Ave | -2.05E03 | **-4.19E03** | -3.24E03 | -3.59E03 | -2.08E03 | -3.25E03 | -3.04E03 | -2.55E03 | -4.06E03 | **-4.19E03** |
|  | Std | 1.04E02 | 2.78E-12 | 2.35E02 | 2.66E02 | 3.84E02 | 2.43E02 | 3.28E02 | 2.02E02 | 6.29E01 | 2.78E-12 |
|  | Best | -2.49E03 | -4.19E03 | -3.83E03 | -3.97E03 | -3.48E03 | -3.71E03 | -3.6E03 | -2.99E03 | -4.19E03 | -4.19E03 |
| F9 | Ave | 6.65E01 | **0** | 1.04E00 | 5.97E-01 | 2.04E01 | 1.67E-01 | 1.7E01 | **0** | 1.24E-02 | **0** |
|  | Std | 3.55E01 | 0 | 1.82E00 | 2.15E00 | 1.02E01 | 9.15E-01 | 1.18E01 | 0 | 6.68E-02 | 0 |
|  | Best | 1.39E01 | 0 | 0 | 0 | 9.95E00 | 0 | 3.98E00 | 0 | 0 | 0 |
| F10 | Ave | 2E01 | 3.76E-15 | 3.88E-15 | **4.44E-16** | 1.46E01 | 3.29E-15 | 1.25E-14 | **4.44E-16** | 2.12E-14 | **4.44E-16** |
|  | Std | 5.01E-03 | 9.01E-16 | 6.49E-16 | 0 | 8.98E00 | 1.45E-15 | 7.95E-15 | 0E00 | 2.32E-14 | 0E00 |
|  | Best | 1.99E01 | 4.44E-16 | 4.44E-16 | 4.44E-16 | 3.98E-07 | 4.44E-16 | 4E-15 | 4.44E-16 | 4E-15 | 4.44E-16 |
| F11 | Ave | 6.38E01 | 3.64E-62 | 1.33E-02 | 3.28E-03 | 2.47E-01 | 1.65E-02 | 1.91E-01 | **0** | 1.32E-03 | 1.58E-52 |
|  | Std | 5.6E)1 | 1.48E-61 | 2.33E-02 | 1.18E-02 | 2.37E-01 | 1.64E-02 | 1.33E-01 | 0 | 2.78E-03 | 8.15E-52 |
|  | Best | 1.09E-19 | 1.03E-72 | 0 | 0 | 6.4E-02 | 1.3E-286 | 3.2E-02 | 0 | 2.15E-08 | 0 |
| F12 | Ave | 2.29E08 | **4.71E-32** | 1.33E-04 | 3.11E-02 | 9.9E-07 | 5E-20 | 1.55E-01 | 4.51E-01 | 5.63E-31 | **4.71E-32** |
|  | Std | 2.73E08 | 1.67E-47 | 5.76E-04 | 1.15E-01 | 3.08E-06 | 1.5E-19 | 3.54E-01 | 1.51E-01 | 2E-30 | 1.67E-47 |
|  | Best | 4.18E-20 | 4.71E-32 | 8.97E-14 | 4.71E-32 | 1.51E-14 | 3.91E-23 | 4.71E-32 | 1.91E-01 | 4.71E-32 | 4.71E-32 |
| F13 | Ave | 3.11E08 | **1.35E-32** | 4.86E-02 | 1.52E-01 | 3.9E-03 | 8.08E-03 | 8.37E-02 | 9.17E-01 | 2.76E-24 | **1.35E-32** |
|  | Std | 3.13E08 | 5.57E-48 | 3.02E-02 | 1.1E-01 | 2.07E-02 | 3.08E-02 | 4.14E-01 | 1.48E-01 | 1.51E-23 | 5.57E-48 |
|  | Best | 2.81E-19 | 1.35E-32 | 3.12E-12 | 1.35E-32 | 4.08E-14 | 2.14E-21 | 1.56E-32 | 4.84E-01 | 2.19E-32 | 1.35E-32 |
| Num Of Tops |  | 0 | 4 | 0 | 1 | 0 | 0 | 0 | 3 | 0 | **5** |
| Rank |  | 5 | 2 | 5 | 4 | 5 | 5 | 5 | 3 | 5 | **1** |

**Table 5** Optimization results for fixed dimension multimodal benchmark functions.

| Fun. | Index | CMA-ES | SHADE | TLBO | LCBO | GMBA | EO | MFO | SO | VLA | AVLA |
|---|---|---|---|---|---|---|---|---|---|---|---|
| F14 | Ave | 2.58E01 | 9.98E-01 | 1.03E00 | 2.37E00 | 9.98E-01 | 3.35E00 | 1.53E00 | 1.61E00 | 9.98E-01 | 9.98E-01 |

|      | Std  | 8.97E01  | 1.13E-16 | 1.81E-01 | 3.37E00  | 1.81E-08 | 3.67E00  | 1.24E00  | 7.83E-01 | 1.13E-16 | 1.13E-16 |
|------|------|----------|----------|----------|----------|----------|----------|----------|----------|----------|----------|
|      | Best | 9.98E-01 | 9.98E-01 | 9.98E-01 | 9.98E-01 | 9.98E-01 | 9.98E-01 | 9.98E-01 | 9.98E-01 | 9.98E-01 | 9.98E-01 |
| F15  | Ave  | 1.44E-02 | 3.07E-04 | 4.42E-04 | 2.34E-03 | 7.06E-04 | 1.68E-03 | 1.24E-03 | 6.72E-04 | 3.51E-04 | 3.09E-04 |
|      | Std  | 8.18E-03 | 1.42E-19 | 3E-04    | 6.11E-03 | 4.13E-04 | 5.08E-03 | 1.43E-03 | 2.36E-04 | 4.61E-05 | 4.07E-06 |
|      | Best | 2.25E-03 | 3.07E-04 | 3.07E-04 | 3.07E-04 | 3.07E-04 | 3.07E-04 | 3.08E-04 | 3.37E-04 | 3.08E-04 | 3.07E-04 |
| F16  | Ave  | 1.39E01  | -1.03E00 | -1.03E00 | -1.03E00 | -1.03E00 | -1.03E00 | -1.03E00 | -1.03E00 | -1.03E00 | -1.03E00 |
|      | Std  | 8.21E01  | 6.78E-16 | 5.86E-16 | 6.71E-16 | 2.95E-11 | 6.78E-16 | 6.25E-16 | 2.31E-04 | 6.78E-16 | 6.78E-16 |
|      | Best | -1.03E00 | -1.03E00 | -1.03E00 | -1.03E00 | -1.03E00 | -1.03E00 | -1.03E00 | -1.03E00 | -1.03E00 | -1.03E00 |
| F17  | Ave  | 1.44E00  | 3.98E-01 | 3.98E-01 | 3.98E-01 | 3.98E-01 | 3.98E-01 | 3.98E-01 | 4.04E-01 | 3.98E-01 | 3.98E-01 |
|      | Std  | 2.12E00  | 0E00     | 0E00     | 0E00     | 1.66E-08 | 0E00     | 0E00     | 7.8E-03  | 0E00     | 0E00     |
|      | Best | 3.98E-01 | 3.98E-01 | 3.98E-01 | 3.98E-01 | 3.98E-01 | 3.98E-01 | 3.98E-01 | 3.98E-01 | 3.98E-01 | 3.98E-01 |
| F18  | Ave  | 4.37E01  | 3E00     | 3E00     | 3E00     | 1.11E01  | 3E00     | 3E00     | 3.03E00  | 3E00     | 3E00     |
|      | Std  | 7.35E01  | 2.23E-15 | 2.26E-15 | 2.44E-15 | 1.26E01  | 3.38E-09 | 3.02E-15 | 3.34E-02 | 2.26E-15 | 2.26E-15 |
|      | Best | 3E00     | 3E00     | 3E00     | 3E00     | 3E00     | 3E00     | 3E00     | 3E00     | 3E00     | 3E00     |
| F19  | Ave  | -3.19E00 | -3.78E00 | -3.78E00 | -3.78E00 | -3.54E00 | -3.78E00 | -3.78E00 | -3.76E00 | -3.78E00 | -3.78E00 |
|      | Std  | 8.24E-01 | 2.26E-15 | 2.26E-15 | 2.26E-15 | 6.17E-01 | 2.26E-15 | 2.26E-15 | 9.2E-03  | 2.26E-15 | 2.26E-15 |
|      | Best | -3.78E00 | -3.78E00 | -3.78E00 | -3.78E00 | -3.78E00 | -3.78E00 | -3.78E00 | -3.78E00 | -3.78E00 | -3.78E00 |
| F20  | Ave  | -5.14E-01| -5.73E-01| -5.73E-01| -5.73E-01| -5.72E-01| -5.73E-01| -5.73E-01| -5.11E-01| -5.73E-01| -5.73E-01|
|      | Std  | 1.11E-01 | 1.92E-16 | 3.39E-16 | 3.57E-16 | 2.63E-04 | 3.22E-12 | 2.14E-16 | 1.45E-02 | 3.43E-16 | 2.4E-16  |
|      | Best | -5.73E-01| -5.73E-01| -5.73E-01| -5.73E-01| -5.73E-01| -5.73E-01| -5.73E-01| -5.59E-01| -5.73E-01| -5.73E-01|
| F21  | Ave  | -4.97E00 | -1.02E01 | -8.86E00 | -6.58E00 | -7.05E00 | -1.02E01 | -6.65E00 | -3.47E00 | -1.02E01 | -1.02E01 |
|      | Std  | 4.17E00  | 8.24E-15 | 2.45E00  | 3.68E00  | 3.27E00  | 5.42E-15 | 3.45E00  | 1.76E00  | 8.42E-15 | 8.76E-15 |
|      | Best | -1.02E01 | -1.02E01 | -1.02E01 | -1.02E01 | -1.02E01 | -1.02E01 | -1.02E01 | -8.12E00 | -1.02E01 | -1.02E01 |
| F22  | Ave  | -4.49E00 | -1.04E01 | -9.6E00  | -7.15E00 | -6.6E00  | -3.7E00  | -8.38E00 | -3.83E00 | -1.04E01 | -1.04E01 |
|      | Std  | 4.15E00  | 3.3E-16  | 1.8E00   | 3.79E00  | 3.5E00   | 3.08E-12 | 3.45E00  | 1.5E00   | 0E00     | 0E00     |
|      | Best | -1.04E01 | -1.04E01 | -1.04E01 | -1.04E01 | -1.04E01 | -3.7E00  | -1.04E01 | -7.61E00 | -1.04E01 | -1.04E01 |
| F23  | Ave  | -3.85E00 | -1.05E01 | -9.79E00 | -7E00    | -4.91E00 | -1.05E01 | -7.61E00 | -3.65E00 | -1.05E01 | -1.05E01 |
|      | Std  | 3.51E00  | 3.57E-15 | 1.98E00  | 3.42E00  | 2.99E00  | 5.42E-15 | 3.45E00  | 1.52E00  | 3.57E-15 | 3.28E-15 |
|      | Best | -1.05E01 | -1.05E01 | -1.05E01 | -1.05E01 | -1.05E01 | -1.05E01 | -1.05E01 | -8.89E00 | -1.05E01 | -1.05E01 |
| Num Of Tops | | 0 | **10** | 5 | 5 | 3 | 7 | 5 | 1 | 9 | 9 |
| Rank |      | 7        | **1**    | 4        | 4        | 5        | 3        | 4        | 6        | 2        | 2        |

**Table 6** Optimization results and comparison for composite benchmark functions.

| Fun. | Index | CMA-ES  | SHADE   | TLBO    | LCBO    | GMBA    | EO      | MFO     | SO      | VLA     | AVLA    |
|------|-------|---------|---------|---------|---------|---------|---------|---------|---------|---------|---------|
| F24  | Ave   | 2.66E02 | 5E01    | 6.33E01 | 5.85E-03| 7.72E-05| 7.34E01 | 3.2E00  | 2.29E02 | 5.34E00 | 6.67E00 |
|      | Std   | 1.48E02 | 5.72E01 | 6.15E01 | 1.93E-02| 3.97E-04| 7.85E01 | 4.39E00 | 5.97E01 | 1.81E01 | 02.54E01|
|      | Best  | 3.5E01  | 0       | 9.35E-11| 0       | 3.19E-15| 3.45E-19| 0E00    | 1.27E02 | 9.63E-29| 00E00   |
| F25  | Ave   | 2.94E02 | 5.63E00 | 5.78E01 | 4.33E01 | 1.03E01 | 2.8E01  | 4.55E01 | 2.58E02 | 1.91E01 | 0       |
|      | Std   | 7.75E01 | 2.14E01 | 4.69E01 | 2.27E01 | 8.23E00 | 4.58E01 | 1.59E01 | 5.12E01 | 1.54E01 | 0       |
|      | Best  | 1.9E02  | 0       | 4.34E00 | 5.38E-29| 2.09E-11| 8.61E-18| 5.39E00 | 1.39E02 | 2.14E-07| 0       |
| F26  | Ave   | 1.68E03 | 3.29E01 | 2.27E02 | 1.57E02 | 4.18E02 | 2.38E02 | 5.18E02 | 9E02    | 3.04E00 | 4E01    |
|      | Std   | 1.55E02 | 9.12E01 | 2.29E02 | 2.24E02 | 2.69E02 | 2.05E02 | 0E00    | 1.43E01 | 1.18E02 |         |
|      | Best  | 1.23E03 |         | 1.78E-05| 0       | 1.07E-12| 1.93E-18| 3.33E-30| 9E02    | 1.7E-20 | 0       |
| F27  | Ave   | 6.32E01 | 1.5E01  | 6.31E00 | 6.31E00 | 6.31E00 | 1.5E01  | 6.31E00 | 8.64E01 | 6.31E00 | 6.31E00 |
|      | Std   | 7.86E01 | 3.31E01 | 3.67E-03| 1.76E-03| 1.84E-05| 3.31E01 | 2.88E-04| 2.11E01 | 2.46E-04| 3.64E-04|
|      | Best  | 6.31E00 | 6.31E00 | 6.31E00 | 6.31E00 | 6.31E00 | 6.31E00 | 6.31E00 | 5.22E01 | 6.31E00 | 6.31E00 |
| F28  | Ave   | 2.69E01 | 3.97E01 | 2.1E01  | 1.9E00  | 5.12E-01| 6.39E01 | 5.44E00 | 6.76E01 | 1.17E00 | 1.85E-01|
|      | Std   | 4.28E01 | 5.71E01 | 4.48E01 | 8.43E-01| 2.86E-01| 6.02E01 | 2.16E01 | 3.47E01 | 3.71E-01| 1.8E-01 |
|      | Best  | 4.24E00 | 0       | 5.66E-01| 7.06E-01| 1.76E-01| 1.46E-01| 6.28E-01| 1.49E01 | 5.05E-01| 0       |
| F29  | Ave   | 1.01E03 | 3.1E02  | 5.29E02 | 6.52E02 | 4.54E02 | 4.74E02 | 7.19E02 | 7.77E02 | 2.32E02 | 3.6E02  |
|      | Best  | 1.07E02 | 2.93E02 | 1.68E02 | 2.91E02 | 3.28E02 | 3.24E02 | 2.2E02  | 5.67E01 | 1.71E02 | 1.92E02 |
|      | Std   | 8.29E02 | 0       | 4E02    | 0       | 5.19E-11| 5.28E00 | 0       | 6.81E02 | 8.82E-18| 0E00    |
| Num Of Tops | | 0 | 0 | 1 | 1 | 1 | 0 | 2 | 0 | **3** | **3** |
| Rank |       | 4       | 4       | 3       | 3       | 3       | 4       | 2       | 4       | **1**   | **1**   |

**Table 7** The final rank for the 29 test problems.

| Fun.        | Index | CMA-ES | SHADE | TLBO | LCBO | GMBA | EO | MFO | SO | VLA | AVLA |
|-------------|-------|--------|-------|------|------|------|----|-----|----|-----|------|
| Num Of Tops |       | 0      | 16    | 9    | 12   | 4    | 8  | 7   | 8  | 13  | **19** |
| Rank        |       | 9      | 2     | 5    | 4    | 8    | 6  | 7   | 6  | 3   | **1** |

The 30 problems in the second group come from the CEC'14 benchmark suite. The dimensions of these problems are all set to 10 in this paper. As well-known, SHADE is the winner of CEC 14. The results given in Table 8 show that among 30 problems, SHADE obtains 23 best solutions for 30 problems. AVLA and VLA obtain 21 and 11 first places with their means, respectively. If we remove VLA from the list of the compared algorithms, AVLA will obtain 23 first places of means. This is a very surprising result. Two observations are worth mentioning here. One is that though AVLA can be viewed as the advanced version of VLA, in some cases VLA can beat AVLA easily, e.g. when solving problems CF4 and CF27 in this group. The other observation is that though the whole performances of SHADE and AVLA are equivalent, they can be distinguished from each other by showing better results for different problems. For example, when dealing with CF9 in this group, SHADE gives the better result. But when dealing with CF27 in this group, the result obtained by AVLA is better than the one given by SHADE.

The third group consists of problems come from the CEC'15 learning based benchmark suite. The dimensions of these problems are all set to 30 in this paper. According to the data given in Table 9, the performances of SHADE, AVLA and VLA are similar to that related to the second group. AVLA lags behind SHADE due to its worse performance when tackling the simple multimodal functions CF3, CF4 and CF5 in this group. This shows the weakness of AVLA when dealing with such type of problem. Just as shown in Table 6, SHADE is weak when applied to the composite benchmark functions in the first group. The above observation reminds us a fair comparison between algorithms should be carried out on various types of problems.

The fourth group includes 12 problems from the CEC'22 test suite. In this paper we only consider the versions of these problems with dimension 20. When used to tackle these 12 problems, SHADE and AVLA share the first place by obtaining 9 first places out of 12 trials. VLA follows them to gain the second rank with 8 best results for 12 problems. We can see VLA outperforms others by giving the best solutions to problems CF2 and CF11 in this group. This demonstrates the special value of VLA as an outstanding meta-heuristic algorithm. The results associated with this group are presented in Table10.

Table 8 Results for optimization problems in CEC 2014.

| Fun. | Index | CMA_ES | SHADE | TLBO | LCBO | GMBA | EO | MFO | SO | VLA | AVLA |
|---|---|---|---|---|---|---|---|---|---|---|---|
| CF1 | Ave | 1.54E08 | 1E02 | 8.17E04 | 4.89E06 | 7.36E04 | 5.34E04 | 2.57E06 | 4.65E07 | 1.8E04 | 1E02 |
| | Std | 9.03E07 | 0E00 | 6.7E04 | 1.88E07 | 3.62E04 | 4.96E04 | 6.44E06 | 2.01E07 | 8.35E03 | 1.08E-10 |
| | Best | 1.42E07 | 1E02 | 4.23E03 | 5.69E02 | 1.3E04 | 1.9E04 | 3.65E04 | 1.24E07 | 3.37E03 | 1E02 |
| CF2 | Ave | 1.3E10 | 2E02 | 9.84E04 | 3.8E02 | 1.08E06 | 2.84E03 | 6.27E03 | 2.8E09 | 4.78E02 | 2E02 |
| | Std | 4.87E09 | 0E00 | 4.1E05 | 2.27E02 | 5.91E06 | 1.02E00 | 4.6E03 | 8.3E08 | 6.76E02 | 0E00 |
| | Best | 2.24E09 | 2E02 | 2.2E02 | 2E02 | 2.24E02 | 2.84E03 | 2.2E02 | 1.4E09 | 2E02 | 2E00 |
| CF3 | Ave | 5.38E06 | 3E02 | 6.25E02 | 6.03E02 | 5.97E04 | 2.63E03 | 2.64E04 | 1.29E04 | 4.07E02 | 3E02 |
| | Std | 1.22E07 | 0E00 | 4.56E02 | 1.06E03 | 1.3E05 | 3.77E03 | 2.28E04 | 2.05E03 | 8.76E01 | 1.49E-14 |
| | Best | 5.58E04 | 3E02 | 3E02 | 3E02 | 4.71E03 | 3E02 | 5.27E02 | 7.16E03 | 3.05E02 | 3E02 |
| CF4 | Ave | 1.5E03 | 4.22E02 | 4.21E02 | 4.23E02 | 4.27E02 | 4.35E02 | 4.36E02 | 1.02E03 | 4.04E02 | 4.15E02 |
| | Std | 5.08E02 | 1.65E01 | 1.93E01 | 1.81E01 | 1.76E01 | 5.63E-02 | 3E01 | 2.48E02 | 8.27E00 | 1.69E01 |
| | Best | 5.9E02 | 4E02 | 4E02 | 4E02 | 4E02 | 4.35E02 | 4E02 | 6.43E02 | 4E02 | 4E02 |
| CF5 | Ave | 5.21E02 | 5.2E02 | 5.2E02 | 5.2E02 | 5.2E02 | 5.2E02 | 5.2E02 | 5.2E02 | 5.2E02 | 5.19E02 |
| | Std | 1.71E-01 | 1.9E00 | 2.44E00 | 7.55E-02 | 2.16E-01 | 2.45E-04 | 1.28E-01 | 7.51E-02 | 2E00 | 3.66E00 |
| | Best | 5.2E02 | 5.2E02 | 5.07E02 | 5.2E02 | 5.2E02 | 5.2E02 | 5.2E02 | 5.2E02 | 5.09E02 | 5E02 |
| CF6 | Ave | 6.11E02 | 6E02 | 6.03E02 | 6.06E02 | 6.05E02 | 6E02 | 6.05E02 | 6.1E02 | 6.01E02 | 6E02 |
| | Std | 1.07E00 | 0E00 | 1.26E00 | 1.59E00 | 1.75E00 | 5.83E-03 | 1.68E00 | 5.51E-01 | 4.25E-01 | 1.93E-01 |
| | Best | 6.08E02 | 6E02 | 6.01E02 | 6.02E02 | 6E02 | 6E02 | 6.02E02 | 6.09E02 | 6E02 | 6E02 |
| CF7 | Ave | 8.24E02 | 7E02 | 7.02E02 | 7.01E02 | 7E02 | 7E02 | 7.03E02 | 7.85E02 | 7E02 | 7E02 |
| | Std | 3.8E01 | 1.15E-02 | 2.21E00 | 1.67E00 | 1.05E-01 | 3.47E-13 | 4.08E00 | 2.17E01 | 1.5E-02 | 1.49E-02 |
| | Best | 7.58E02 | 7E02 | 7E02 | 7E02 | 7E02 | 7E02 | 7E02 | 7.4E02 | 7E02 | 7E02 |
| CF8 | Ave | 8.78E02 | 8E02 | 8.16E02 | 8.26E02 | 8.21E02 | 8.08E02 | 8.19E02 | 8.68E02 | 8E02 | 8E02 |
| | Std | 2.31E01 | 0E00 | 6.7E00 | 7.93E00 | 1.61E01 | 1.79E00 | 7.95E00 | 4.61E00 | 4.43E-01 | 0E00 |
| | Best | 8.42E02 | 8E02 | 8.04E02 | 8.1E02 | 8.06E02 | 8.08E02 | 8.07E02 | 8.57E02 | 8E02 | 8E02 |
| CF9 | Ave | 9.96E02 | 9.02E02 | 9.21E02 | 9.28E02 | 9.26E02 | 9.34E02 | 9.66E02 | 9.07E02 | 9.04E02 | |
| | Std | 2.35E01 | 9.96E-01 | 4.83E00 | 8.68E00 | 1.14E01 | 6.01E00 | 1.2E01 | 5.79E00 | 1.35E00 | 1.66E00 |
| | Best | 9.25E02 | 9.01E02 | 9.11E02 | 9.16E02 | 9.09E02 | 9.05E02 | 9.1E02 | 9.49E02 | 9.05E02 | 9.01E02 |
| CF10 | Ave | 2.35E03 | 1E03 | 1.34E03 | 1.48E03 | 1.21E03 | 1.29E03 | 1.43E03 | 2.48E03 | 1.02E03 | 1E03 |
| | Std | 1.98E02 | 5.1E-02 | 2.5E02 | 2.46E02 | 2.41E02 | 9.5E01 | 2.56E02 | 1.36E02 | 1.85E01 | 3.16E00 |
| | Best | 2.01E03 | 1E03 | 1.01E03 | 1.06E03 | 1.05E03 | 1.15E03 | 1.03E03 | 2.21E03 | 1.01E03 | 1E03 |
| CF11 | Ave | 3.27E03 | 1.18E03 | 1.94E03 | 1.9E03 | 1.86E03 | 1.66E03 | 2.17E03 | 2.61E03 | 1.35E03 | 1.29E03 |
| | Std | 1.91E02 | 1E02 | 2.52E02 | 2.93E02 | 2.26E02 | 2.66E02 | 3.77E02 | 1.65E02 | 7.6E01 | 1.15E02 |
| | Best | 2.76E03 | 1.1E03 | 1.37E03 | 1.35E03 | 1.48E03 | 1.22E03 | 1.54E03 | 2.05E03 | 1.22E03 | 1.11E03 |
| CF12 | Ave | 1.2E03 | 1.2E03 | 1.2E03 | 1.2E03 | 1.2E03 | 1.2E03 | 1.2E03 | 1.2E03 | 1.2E03 | 1.2E03 |
| | Std | 5.92E-01 | 9.29E-02 | 2.97E-01 | 3.45E-01 | 3.66E-01 | 1.48E-02 | 2.1E-01 | 2.95E-01 | 5.1E-02 | 9.15E-02 |
| | Best | 1.2E03 | 1.2E03 | 1.2E03 | 1.2E03 | 1.2E03 | 1.2E03 | 1.2E03 | 1.2E03 | 1.2E03 | 1.2E03 |
| CF13 | Ave | 1.3E03 | 1.3E03 | 1.3E03 | 1.3E03 | 1.3E03 | 1.3E03 | 1.3E03 | 1.3E03 | 1.3E03 | 1.3E03 |
| | Std | 1.23E00 | 2.9E-02 | 1.42E-01 | 1.56E-01 | 2.16E-01 | 4.63E-13 | 1.2E-01 | 4.56E-01 | 3.04E-01 | 2.46E-02 |
| | Best | 1.3E03 | 1.3E03 | 1.3E03 | 1.3E03 | 1.3E03 | 1.3E03 | 1.3E03 | 1.3E03 | 1.3E03 | 1.3E03 |
| CF14 | Ave | 1.44E03 | 1.4E03 | 1.4E03 | 1.4E03 | 1.4E03 | 1.4E03 | 1.4E03 | 1.42E03 | 1.4E03 | 1.4E03 |
| | Std | 1.66E01 | 3.85E-02 | 8.14E-01 | 2.64E-01 | 2.37E-01 | 4.18E-02 | 3.52E-01 | 4.44E00 | 4.78E-02 | 5.27E-02 |
| | Best | 1.41E03 | 1.4E03 | 1.4E03 | 1.4E03 | 1.4E03 | 1.4E03 | 1.4E03 | 1.41E03 | 1.4E03 | 1.4E03 |
| CF15 | Ave | 2.08E06 | 1.5E03 | 1.5E03 | 1.5E03 | 1.5E03 | 1.5E03 | 1.52E03 | 1.62E03 | 1.5E03 | 1.5E03 |
| | Std | 2.31E06 | 1.44E-01 | 2.18E00 | 2.38E00 | 1.23E00 | 1.02E-02 | 4.01E01 | 8.15E01 | 1.26E-01 | 1.57E-01 |
| | Best | 2.63E03 | 1.5E03 | 1.5E03 | 1.5E03 | 1.5E03 | 1.5E03 | 1.5E03 | 1.52E03 | 1.5E03 | 1.5E03 |
| CF16 | Ave | 1.6E03 | 1.6E03 | 1.6E03 | 1.6E03 | 1.6E03 | 1.6E03 | 1.6E03 | 1.6E03 | 1.6E03 | 1.6E03 |
| | Std | 1.57E-01 | 5.51E-01 | 3.52E-01 | 3.96E-01 | 6.04E-01 | 2.31E-04 | 2.86E-01 | 1.74E-01 | 2.57E-01 | 5.39E-01 |
| | Best | 1.6E03 | 1.6E03 | 1.6E03 | 1.6E03 | 1.6E03 | 1.6E03 | 1.6E03 | 1.6E03 | 1.6E03 | 1.6E03 |
| CF17 | Ave | 9.89E06 | 1.71E03 | 4.11E03 | 5.77E03 | 3.84E05 | 3.82E03 | 3.42E04 | 2.52E05 | 5.18E03 | 1.73E03 |
| | Std | 1.49E07 | 4.05E01 | 2E03 | 3.93E03 | 6.06E05 | 1.21E05 | 7.3E04 | 1.21E05 | 2.37E03 | 4.14E01 |
| | Best | 1.32E05 | 1.7E03 | 1.85E03 | 1.9E03 | 2.3E03 | 3.75E03 | 2.23E03 | 3.3E04 | 2.25E03 | 1.7E03 |
| CF18 | Ave | 1.29E08 | 1.8E03 | 4.69E03 | 8.37E03 | 2.95E04 | 2.89E03 | 1.53E04 | 3.65E04 | 2.6E03 | 1.8E03 |
| | Std | 2.56E08 | 1.7E-01 | 2.61E03 | 5.9E03 | 9.48E04 | 1.16E03 | 1.37E04 | 1.91E04 | 8.3E02 | 1.17E00 |
| | Best | 3.09E03 | 1.8E03 | 1.82E03 | 2.3E03 | 1.92E03 | 1.93E03 | 2.13E03 | 1.26E04 | 1.86E03 | 1.8E03 |

| | | | | | | | | | | |
|---|---|---|---|---|---|---|---|---|---|---|
| CF19 | Ave | 1.91E03 | 1.9E03 | 1.9E03 | 1.9E03 | 1.9E03 | 1.9E03 | 1.9E03 | 1.91E03 | 1.9E03 | 1.9E03 |
| | Std | 6.25E00 | 2.87E-02 | 6.37E-01 | 1.73E00 | 1.34E00 | 2.37E-02 | 1.12E00 | 4.2E00 | 2.77E-01 | 2.39E-01 |
| | Best | 1.91E03 | 1.9E03 | 1.9E03 | 1.9E03 | 1.9E03 | 1.9E03 | 1.9E03 | 1.91E03 | 1.9E03 | 1.9E03 |
| CF20 | Ave | 2.22E06 | 2E03 | 2.12E03 | 2.3E03 | 3.56E03 | 2.14E03 | 1.54E04 | 9.25E03 | 2.3E03 | 2E03 |
| | Std | 1.11E07 | 5.71E-02 | 8.29E01 | 5.98E02 | 2.64E03 | 5.27E02 | 1.36E04 | 3.02E03 | 2.94E02 | 2.49E-01 |
| | Best | 4.88E03 | 2E03 | 2.01E03 | 2.02E03 | 2.06E03 | 2.04E03 | 2.01E03 | 3.24E03 | 2.01E03 | 2E03 |
| CF21 | Ave | 3.34E06 | 2.1E03 | 2.21E03 | 4.59E03 | 1.09E06 | 2.41E03 | 8.25E03 | 4.25E04 | 2.67E03 | 2.1E03 |
| | Std | 5.13E05 | 2.49E-01 | 9.66E01 | 7.46E03 | 1.42E06 | 9.23E00 | 8.05E03 | 2.71E04 | 4.42E02 | 2.19E01 |
| | Best | 1.47E06 | 2.1E03 | 2.11E03 | 2.13E03 | 6.7E03 | 2.4E03 | 2.45E03 | 5.54E03 | 2.15E03 | 2.1E03 |
| CF22 | Ave | 2.62E03 | 2.2E03 | 2.24E03 | 2.29E03 | 2.24E03 | 2.28E03 | 2.26E03 | 2.39E03 | 2.21E03 | 2.2E03 |
| | Std | 1.86E02 | 9.34E-02 | 2.4E01 | 7.17E01 | 3.97E01 | 5.47E01 | 5.76E01 | 3.55E01 | 7.76E00 | 5.02E00 |
| | Best | 2.29E03 | 2.2E03 | 2.22E03 | 2.22E03 | 2.22E03 | 2.22E03 | 2.22E03 | 2.31E03 | 2.2E03 | 2.2E03 |
| CF23 | Ave | 2.97E03 | 2.63E03 | 2.63E03 | 2.54E03 | 2.63E03 | 2.63E03 | 2.64E03 | 2.5E03 | 2.63E03 | 2.57E03 |
| | Std | 1.84E02 | 9.25E-13 | 1.96E-01 | 6.33E01 | 8.79E-12 | 6.65E-05 | 8.58E00 | 0E00 | 1.62E-10 | 6.52E01 |
| | Best | 2.66E03 | 2.63E03 | 2.63E03 | 2.5E03 | 2.63E03 | 2.63E03 | 2.63E03 | 2.5E03 | 2.63E03 | 2.5E03 |
| CF24 | Ave | 2.64E03 | 2.51E03 | 2.54E03 | 2.59E03 | 2.55E03 | 2.54E03 | 2.54E03 | 2.59E03 | 2.52E03 | 2.51E03 |
| | Std | 3.23E01 | 3.44E00 | 2.71E01 | 2.6E01 | 3.27E01 | 2.81E01 | 1.22E01 | 7.53E00 | 2.95E00 | 2.58E00 |
| | Best | 2.58E03 | 2.5E03 | 2.52E03 | 2.52E03 | 2.57E03 | 2.52E03 | 2.52E03 | 2.57E03 | 2.51E03 | 2.51E03 |
| CF25 | Ave | 2.72E03 | 2.64E03 | 2.66E03 | 2.68E03 | 2.7E03 | 2.7E03 | 2.7E03 | 2.7E03 | 2.63E03 | 2.62E03 |
| | Std | 5.9E00 | 4E01 | 2.34E01 | 1.37E01 | 1.89E00 | 3.04E-04 | 1.47E01 | 2.74E00 | 7.7E00 | 5.85E00 |
| | Best | 2.71E03 | 2.6E03 | 2.62E03 | 2.66E03 | 2.7E03 | 2.7E03 | 2.64E03 | 2.69E03 | 2.62E03 | 2.6E03 |
| CF26 | Ave | 2.7E03 | 2.7E03 | 2.7E03 | 2.7E03 | 2.7E03 | 2.7E03 | 2.7E03 | 2.7E03 | 2.7E03 | 2.7E03 |
| | Std | 1.59E00 | 3.29E-02 | 1.24E-01 | 1.76E-01 | 2.77E-01 | 1.39E-12 | 1.45E-01 | 6.26E-01 | 4.13E-02 | 5E-02 |
| | Best | 2.7E03 | 2.7E03 | 2.7E03 | 2.7E03 | 2.7E03 | 2.7E03 | 2.7E03 | 2.7E03 | 2.7E03 | 2.7E03 |
| CF27 | Ave | 3.32E03 | 2.87E03 | 2.91E03 | 3.08E03 | 3.1E03 | 2.9E03 | 3.08E03 | 2.9E03 | 2.71E03 | 2.72E03 |
| | Std | 4.05E01 | 1.75E02 | 1.97E02 | 1.34E02 | 8.07E01 | 1.66E02 | 1.29E02 | 1.07E02 | 1.94E01 | 7.55E01 |
| | Best | 3.25E03 | 2.7E03 | 2.7E03 | 2.71E03 | 2.7E03 | 2.7E03 | 2.7E03 | 2.85E03 | 2.7E03 | 2.7E03 |
| CF28 | Ave | 3.28E03 | 3.2E03 | 3.25E03 | 3.34E03 | 3.25E03 | 3.29E03 | 3.19E03 | 3E03 | 3.18E03 | 3.18E03 |
| | Std | 8.21E01 | 4.92E01 | 5.63E01 | 1.47E02 | 7.44E01 | 6.81E-07 | 1.12E01 | 0E00 | 6.67E00 | 1.29E01 |
| | Best | 3.23E03 | 3.16E03 | 3.17E03 | 2.9E03 | 3.17E03 | 3.29E03 | 3.17E03 | 3E03 | 3.16E03 | 3.16E03 |
| CF29 | Ave | 1.15E06 | 3.12E03 | 8.88E03 | 2.8E05 | 1.79E05 | 3.35E03 | 3.96E03 | 3.1E03 | 3.2E03 | 3.16E03 |
| | Std | 7.68E05 | 8.37E-01 | 1.93E04 | 7.18E05 | 5.25E05 | 1.01E01 | 5.44E02 | 0E00 | 3.49E01 | 6.31E01 |
| | Best | 4.04E03 | 3.12E03 | 3.16E03 | 3.17E03 | 3.37E03 | 3.34E03 | 3.32E03 | 3.1E03 | 3.14E03 | 3.12E03 |
| CF30 | Ave | 1.46E04 | 3.48E03 | 3.98E03 | 4.25E03 | 4.4E03 | 4.23E03 | 3.89E03 | 4.06E03 | 3.56E03 | 3.51E03 |
| | Best | 1.5E04 | 1.99E01 | 3.74E02 | 4.57E02 | 7.34E02 | 3.5E-01 | 3.49E02 | 2.19E03 | 5.64E01 | 1.33E01 |
| | Std | 4.58E03 | 3.45E03 | 3.56E03 | 3.66E03 | 3.49E03 | 4.23E03 | 3.48E03 | 3.2E03 | 3.48E03 | 3.47E03 |
| Num Of Tops | | 4 | 23 | 7 | 7 | 8 | 9 | 6 | 7 | 11 | 21 |
| Rank | | 8 | 1 | 6 | 6 | 5 | 4 | 7 | 6 | 3 | 2 |

**Table 9** Results for optimization problems in CEC 2015 learning-based test suite.

| Fun. | Index | CMA_ES | SHADE | TLBO | LCBO | GMBA | EO | MFO | SO | VLA | AVLA |
|---|---|---|---|---|---|---|---|---|---|---|---|
| CF1 | Ave | 2.84E07 | 1E02 | 1.6E03 | 1.09E04 | 7.6E03 | 1.11E04 | 3.54E04 | 3.13E04 | 3.63E02 | 1E02 |
| | Std | 1.3E08 | 0E00 | 2.04E03 | 1.77E04 | 2.73E05 | 7.5E03 | 2.59E04 | 5.44E03 | 2.75E02 | 2.64E-15 |
| | Best | 2.52E04 | 1E02 | 1.24E02 | 1.09E02 | 1.01E02 | 4.31E02 | 2.75E03 | 2.22E04 | 1.05E02 | 1E02 |
| CF2 | Ave | 1.42E10 | 2E02 | 3.39E06 | 1.3E04 | 3.75E07 | 8.2E05 | 1.64E08 | 5.8E09 | 3.11E03 | 2E02 |
| | Std | 6.47E09 | 0E00 | 1.83E07 | 1.38E04 | 2.05E08 | 4.83E08 | 1.11E06 | 4.83E08 | 1.81E09 | 1.66E02 | 0E00 |
| | Best | 5.59E09 | 2E02 | 1.07E03 | 2.13E02 | 2.33E02 | 6.52E02 | 6.9E02 | 1.98E09 | 4.35E02 | 2E02 |
| CF3 | Ave | 3.21E02 | 3.19E02 | 3.2E02 | 3.2E02 | 3.2E02 | 3.2E02 | 3.2E02 | 3.2E02 | 3.2E02 | 3.2E02 |
| | Std | 1.4E-01 | 3.28E00 | 8.87E-02 | 5.03E-02 | 1.68E-01 | 1.56E-01 | 1.09E-01 | 7.45E-02 | 6.56E-01 | 1.81E-02 |
| | Best | 3.2E02 | 3.08E02 | 3.2E02 | 3.2E02 | 3.2E02 | 3.2E02 | 3.2E02 | 3.2E02 | 3.17E02 | 3.2E02 |
| CF4 | Ave | 4.83E02 | 4.01E02 | 4.24E02 | 4.28E02 | 4.27E02 | 4.17E02 | 4.32E02 | 4.7E02 | 4.08E02 | 4.04E02 |
| | Std | 3.87E01 | 1E00 | 7.62E00 | 1.12E01 | 1.64E01 | 6.81E00 | 1.41E01 | 8.05E00 | 1.59E00 | 1.6E00 |
| | Best | 4.12E02 | 4E02 | 4.11E02 | 4.08E02 | 4.09E02 | 4.05E02 | 4.1E02 | 4.51E02 | 4.05E02 | 4.01E02 |
| CF5 | Ave | 2.74E03 | 5.42E02 | 1.11E03 | 1.29E03 | 1.29E03 | 1.01E03 | 1.42E03 | 2.04E03 | 7.17E02 | 6.32E02 |
| | Std | 2.03E02 | 4.76E01 | 2.59E02 | 2.5E02 | 4.1E02 | 2.46E02 | 2.7E02 | 1.21E02 | 6.55E01 | 9.65E01 |
| | Best | 2.33E03 | 5E02 | 6.59E02 | 8.57E02 | 8.66E02 | 5.13E02 | 9.43E02 | 1.73E03 | 5.84E02 | 5.07E02 |
| CF6 | Ave | 9.07E05 | 6.05E02 | 1.17E03 | 5.65E03 | 4.15E04 | 3.51E03 | 6.62E03 | 9.6E03 | 1.25E03 | 6.13E02 |
| | Std | 3.63E06 | 2.19E01 | 2.46E02 | 7.01E03 | 4.99E04 | 2.34E03 | 4.8E03 | 4.72E03 | 2.43E02 | 3.29E01 |
| | Best | 1.28E04 | 6E02 | 7.81E02 | 8.42E02 | 4.37E03 | 9.98E02 | 1.31E03 | 2.34E03 | 7.36E02 | 6E02 |
| CF7 | Ave | 7.5E02 | 7E02 | 7.02E02 | 7.03E02 | 7.04E02 | 7.02E02 | 7.03E02 | 7.12E02 | 7.01E02 | 7E02 |
| | Std | 3.79E01 | 8.68E-02 | 8.48E-01 | 1.42E00 | 1.97E00 | 8.44E-01 | 1.31E00 | 2.44E00 | 2.59E-01 | 4.53E-02 |
| | Best | 7.13E02 | 7E02 | 7E02 | 7.01E02 | 7.01E02 | 7.01E02 | 7.01E02 | 7.07E02 | 7E02 | 7E02 |
| CF8 | Ave | 5.58E05 | 8E02 | 1.16E03 | 2.86E03 | 1.93E03 | 1.67E04 | 1.93E03 | 9.86E03 | 3.54E03 | 1.07E03 | 8.02E02 |
| | Std | 1.62E06 | 2.19E-01 | 1.92E02 | 4.84E03 | 5.74E04 | 1.25E03 | 8.41E03 | 1E03 | 1.86E02 | 1.06E01 |
| | Best | 8.91E03 | 8E02 | 8.83E02 | 9.65E02 | 8.27E02 | 8.49E02 | 1.53E03 | 1.97E03 | 8.02E02 | 8E02 |
| CF9 | Ave | 1.04E03 | 1E03 | 1E03 | 1E03 | 1E03 | 1E03 | 1E03 | 1.03E03 | 1E03 | 1E03 |
| | Std | 2.2E01 | 3.74E-02 | 1.76E-01 | 3.93E-01 | 3.28E-01 | 1.33E-01 | 7.28E-01 | 8.8E00 | 4.15E-02 | 5.83E-02 |
| | Best | 1.01E03 | 1E03 | 1E03 | 1E03 | 1E03 | 1E03 | 1E03 | 1.01E03 | 1E03 | 1E03 |
| CF10 | Ave | 3.8E05 | 1.22E03 | 1.46E03 | 2.4E03 | 2.11E04 | 1.81E03 | 7.48E03 | 8.26E03 | 1.4E03 | 1.22E03 |
| | Std | 8.92E05 | 5.12E00 | 1.41E02 | 2.64E03 | 9.35E04 | 4.11E02 | 5.26E03 | 3.37E03 | 7.7E01 | 1.13E00 |
| | Best | 1.06E04 | 1.22E03 | 1.28E03 | 1.28E03 | 1.36E03 | 1.3E03 | 1.44E03 | 4.55E03 | 1.28E03 | 1.22E03 |
| CF11 | Ave | 1.44E03 | 1.29E03 | 1.36E03 | 1.39E03 | 1.4E03 | 1.38E03 | 1.38E03 | 1.34E03 | 1.11E03 | 1.12E03 |
| | Std | 7.99E01 | 1.46E02 | 1.02E02 | 5.36E01 | 7.74E-02 | 7.51E01 | 7.44E01 | 5.79E01 | 2.27E00 | 7.52E01 |
| | Best | 1.4E03 | 1.1E03 | 1.11E03 | 1.11E03 | 1.4E03 | 1.1E03 | 1.11E03 | 1.2E03 | 1.11E03 | 1.1E03 |
| CF12 | Ave | 1.36E03 | 1.3E03 | 1.31E03 | 1.31E03 | 1.31E03 | 1.3E03 | 1.31E03 | 1.33E03 | 1.3E03 | 1.3E03 |
| | Std | 2.17E01 | 3.66E-01 | 2.03E00 | 2.02E00 | 1.8E01 | 9.62E-01 | 2.88E00 | 4.78E00 | 3.76E-01 | 6.23E-01 |
| | Best | 1.32E03 | 1.3E03 | 1.3E03 | 1.3E03 | 1.3E03 | 1.3E03 | 1.3E03 | 1.32E03 | 1.3E03 | 1.3E03 |
| CF13 | Ave | 1.35E03 | 1.33E03 | 1.34E03 | 1.34E03 | 1.34E03 | 1.33E03 | 1.34E03 | 1.35E03 | 1.33E03 | 1.33E03 |
| | Std | 2.01E00 | 1.89E00 | 3.32E00 | 4.25E00 | 4.63E00 | 3.47E00 | 3.1E00 | 2.35E00 | 1.3E00 | 1.59E00 |
| | Best | 1.34E03 | 1.32E03 | 1.33E03 | 1.33E03 | 1.33E03 | 1.33E03 | 1.33E03 | 1.35E03 | 1.32E03 | 1.32E03 |
| CF14 | Ave | 1.1E04 | 5.53E03 | 5.48E03 | 6.85E03 | 7.3E03 | 7.69E03 | 7.67E03 | 1.19E04 | 2.12E03 | 5.07E03 |
| | Std | 3.51E03 | 2.09E03 | 2.44E03 | 2.9E03 | 1.71E03 | 2.27E03 | 1.55E03 | 1.49E03 | 9.78E02 | 2.06E03 |
| | Best | 8.04E03 | 1.5E03 | 1.5E03 | 4.25E03 | 4.25E03 | 1.5E03 | 4.25E03 | 9.28E03 | 1.5E03 | 1.5E03 |
| CF15 | Ave | 1.1E04 | 1.6E03 | 1.6E03 | 1.6E03 | 1.6E03 | 1.6E03 | 1.61E03 | 1.72E03 | 1.6E03 | 1.6E03 |

|  | Std | 1.28E04 | 0E00 | 6.08E00 | 1.13E00 | 4.47E-01 | 2.66E00 | 1.6E01 | 4.01E01 | 8.97E-13 | 4.22E-14 |
|  | Best | 2.9E03 | 1.6E03 | 1.6E03 | 1.6E03 | 1.6E03 | 1.6E03 | 1.6E03 | 1.67E03 | 1.6E03 | 1.6E03 |
| Num Of Tops |  | 0 | 13 | 2 | 2 | 2 | 4 | 1 | 0 | 6 | 8 |
| Rank |  | 7 | 1 | 5 | 5 | 5 | 4 | 6 | 7 | 3 | 2 |

**Table 10** Results for optimization problems in CEC 2022.

| Fun. | Index | CMA_ES | SHADE | TLBO | LCBO | GMBA | EO | MFO | SO | VLA | AVLA |
|---|---|---|---|---|---|---|---|---|---|---|---|
| CF1 | Ave | 5.94E04 | 3E02 | 3E02 | 3E02 | 3E02 | 3.12E02 | 1.77E03 | 5.1E03 | 3E02 | 3E02 |
|  | Std | 1.4E04 | 0E00 | 4.97E-01 | 3E-10 | 4.63E-01 | 2.43E01 | 3.63E03 | 9.54E02 | 1.43E-13 | 0E00 |
|  | Best | 3.52E04 | 3E02 | 3E02 | 3E02 | 3E02 | 3E02 | 3E02 | 2.88E03 | 3E02 | 3E02 |
| CF2 | Ave | 9E02 | 4.03E02 | 4.09E02 | 4.22E02 | 4.11E02 | 4.14E02 | 4.14E02 | 9.08E02 | 4E02 | 4.04E02 |
|  | Std | 2.11E02 | 3.13E00 | 1.93E01 | 2.87E01 | 1.5E01 | 1.89E01 | 1.36E01 | 1.8E02 | 5.12E-02 | 3.77E00 |
|  | Best | 4.42E02 | 4E02 | 4E02 | 4E02 | 4E02 | 4E02 | 4.04E02 | 6.04E02 | 4E02 | 4E02 |
| CF3 | Ave | 6.04E02 | 6.02E02 | 6.03E02 | 6.03E02 | 6.03E02 | 6.02E02 | 6.03E02 | 6.04E02 | 6.03E02 | 6.02E02 |
|  | Std | 2.26E-01 | 4.96E-01 | 3.51E-01 | 4.77E-01 | 5.64E-01 | 5.04E-01 | 4.01E-01 | 1.12E-01 | 1.9E-01 | 5.12E-01 |
|  | Best | 6.04E02 | 6.01E02 | 6.02E02 | 6.02E02 | 6.02E02 | 6.02E02 | 6.03E02 | 6.04E02 | 6.02E02 | 6.01E02 |
| CF4 | Ave | 8.87E02 | 8.03E02 | 8.21E02 | 8.2E02 | 8.33E02 | 8.15E02 | 8.3E02 | 8.51E02 | 8.09E02 | 8.05E02 |
|  | Std | 3.45E01 | 1.32E00 | 8.13E00 | 6.34E00 | 1.78E01 | 5.94E00 | 1.35E01 | 5.89E00 | 2.04E00 | 1.63E00 |
|  | Best | 8.22E02 | 8E02 | 8.09E02 | 8.08E02 | 8.08E02 | 8.05E02 | 8.03E02 | 8.37E02 | 8.04E02 | 8.02E02 |
| CF5 | Ave | 9.15E02 | 9E02 | 9.01E02 | 9.01E02 | 9.03E02 | 9E02 | 9.01E02 | 9.03E02 | 9E02 | 9E02 |
|  | Std | 5.2E00 | 0E00 | 6E-01 | 6.2E-01 | 2.98E00 | 4.02E-01 | 1.11E00 | 4.49E-01 | 2.55E-02 | 0E00 |
|  | Best | 9.08E02 | 9E02 | 9E02 | 9.01E02 | 9E02 | 9E02 | 9E02 | 9.02E02 | 9E02 | 9E02 |
| CF6 | Ave | 6.86E08 | 1.8E03 | 3.6E03 | 4.57E03 | 4.34E03 | 4.59E03 | 5.23E03 | 3.21E07 | 1.91E03 | 1.8E03 |
|  | Std | 4.9E08 | 2.8E-01 | 1.78E03 | 2.28E03 | 1.89E03 | 2.11E03 | 2.46E03 | 2.14E07 | 1.19E02 | 4.22E-01 |
|  | Best | 5.31E03 | 1.8E03 | 1.81E03 | 1.93E03 | 2.03E03 | 1.91E03 | 1.85E03 | 4.26E06 | 1.83E03 | 1.8E03 |
| CF7 | Ave | 2.15E03 | 2E03 | 2.02E03 | 2.03E03 | 2.02E03 | 2.02E03 | 2.02E03 | 2.05E03 | 2E03 | 2E03 |
|  | Std | 4.72E01 | 3.13E-03 | 7.08E00 | 1.05E01 | 6.4E00 | 6.35E00 | 3.58E00 | 6.79E00 | 1.06E00 | 4.28E-01 |
|  | Best | 2.07E03 | 2E03 | 2E03 | 2E03 | 2E03 | 2E03 | 2.02E03 | 2.04E03 | 2E03 | 2E03 |
| CF8 | Ave | 2.31E03 | 2.2E03 | 2.22E03 | 2.22E03 | 2.23E03 | 2.22E03 | 2.22E03 | 2.25E03 | 2.21E03 | 2.2E03 |
|  | Std | 7.81E01 | 1.64E00 | 5.58E00 | 4.62E00 | 1.22E01 | 3.63E00 | 5.06E00 | 7.84E00 | 5.42E00 | 6.01E00 |
|  | Best | 2.23E03 | 2.2E03 | 2.2E03 | 2.2E03 | 2.21E03 | 2.21E03 | 2.22E03 | 2.24E03 | 2.2E03 | 2.2E03 |
| CF9 | Ave | 3E03 | 2.53E03 | 2.53E03 | 2.54E03 | 2.53E03 | 2.53E03 | 2.53E03 | 2.71E03 | 2.53E03 | 2.53E03 |
|  | Std | 3.72E02 | 0E00 | 9.04E-02 | 2.04E01 | 2.68E01 | 1.81E-10 | 5.99E00 | 3.34E01 | 1.11E-11 | 0E00 |
|  | Best | 2.53E03 | 2.53E03 | 2.53E03 | 2.53E03 | 2.53E03 | 2.53E03 | 2.53E03 | 2.66E03 | 2.53E03 | 2.53E03 |
| CF10 | Ave | 3.94E03 | 2.51E03 | 2.51E03 | 2.58E03 | 2.56E03 | 2.54E03 | 2.52E03 | 2.55E03 | 2.5E03 | 2.5E03 |
|  | Std | 2.91E02 | 3.23E01 | 3.67E01 | 6.41E01 | 1.72E02 | 5.63E01 | 4.46E01 | 2.26E01 | 6.34E-02 | 6.83E-02 |
|  | Best | 2.72E03 | 2.5E03 | 2.5E03 | 2.5E03 | 2.5E03 | 2.5E03 | 2.5E03 | 2.51E03 | 2.5E03 | 2.5E03 |
| CF11 | Ave | 4.74E03 | 3.74E03 | 3.87E03 | 4.23E03 | 4.33E03 | 4.2E03 | 4.49E03 | 5.16E03 | 3.03E03 | 3.74E03 |
|  | Std | 3.16E01 | 8.44E02 | 8.04E02 | 7.16E02 | 6.84E02 | 7.61E02 | 5.37E02 | 5.54E02 | 2.78E01 | 8.44E02 |
|  | Best | 4.7E03 | 3E03 | 3.02E03 | 3.03E03 | 3.01E03 | 3.02E03 | 3.08E03 | 3.82E03 | 2.98E03 | 3E03 |
| CF12 | Ave | 289E03 | 2.86E03 | 2.87E03 | 2.88E03 | 2.89E03 | 2.86E03 | 2.86E03 | 2.98E03 | 2.86E03 | 2.86E03 |
|  | Std | 3.11E01 | 1.34E00 | 4.39E00 | 2.36E00 | 1.53E00 | 1.29E00 | 1.3E00 | 3.84E01 | 1.27E00 | 1.23E00 |
|  | Best | 2.87E03 | 2.86E03 | 2.86E03 | 2.86E03 | 2.86E03 | 2.86E03 | 2.86E03 | 2.91E03 | 2.86E03 | 2.86E03 |
| Num Of Tops |  | 0 | 9 | 2 | 1 | 2 | 4 | 2 | 0 | 8 | 9 |
| Rank |  | 6 | 1 | 4 | 5 | 4 | 3 | 4 | 6 | 2 | 1 |

The last group is composed of 14 well-known engineering optimization problems. In Table 11, "EP" stands for engineering problem. The problems indicated by EP1 to EP14 are the gear chain design problem, the three-bar truss design, the robot gripper design, the Belleville spring design, the step-cone pulley design, the hydrostatic thrust bearing design, the car crashworthiness problem, the cantilever beam design, the rolling element bearing design, the multiple disc clutch brake design, the speed reducer design, the tension/compression spring design, the welded beam design, and the pressure vessel design, respectively. In contrast to the known optima of CEC test suite, the optimal results of these engineering optimization problems are unknown and harder to be predicted. All these 14 problems have been studied by many researchers in the past. It is easy for readers to gain the related information of these problems by Googling them or checking related papers such as MLA [64]. When applied to these engineering problems, AVLA outperforms others and obtains 13 first places of the means out of 14. SHADE gains the second rank by providing 11 best means of the results for 15 problems. VLA follows SHADE and supplies 10 best means of the results. Among the left algorithms, TLBO occupies the fourth rank by obtaining 6 best means.

**Table 11** The optimal results for 14 engineering optimization problems.

| Fun. | Index | CMA_ES | SHADE | TLBO | LCBO | GMBA | EO | MFO | SO | VLA | AVLA |
|---|---|---|---|---|---|---|---|---|---|---|---|
| EP1 | Ave | 5.87E-05 | 7.46E-12 | 4.1E-10 | 3.45E-08 | 4.86E-08 | 4.39E-10 | 1.55E-08 | 2.31E-07 | 8.81E-12 | 4.8E-10 |
|  | Std | 1.31E-04 | 8.77E-12 | 4.65E-10 | 9.51E-08 | 2.52E-07 | 4.79E-10 | 1.79E-08 | 3.52E-07 | 9.5E-12 | 5.04E-10 |
|  | Best | 2.36E-09 | 2.7E-12 | 2.7E-12 | 3.07E-10 | 2.7E-12 | 2.7E-12 | 1.26E-09 | 2.36E-09 | 2.7E-12 | 2.7E-12 |
| EP2 | Ave | 2.82E02 | 2.64E02 | 2.64E02 | 2.64E02 | 2.71E02 | 2.64E02 | 2.64E02 | 2.64E02 | 2.64E02 | 2.64E02 |
|  | Std | 3.22E00 | 2.59E-14 | 1.23E-05 | 1.44E-05 | 6.7E00 | 3.13E-04 | 1.65E-01 | 4.43E-01 | 2.86E-05 | 1.89E-12 |
|  | Best | 2.69E02 | 2.64E02 | 2.64E02 | 2.64E02 | 2.64E02 | 2.64E02 | 2.64E02 | 2.64E02 | 2.64E02 | 2.64E02 |
| EP3 | Ave | 3.77E03 | 4.32E00 | 4.94E00 | 6.92E00 | 5.54E00 | 4.5E00 | 5.68E00 | 8.66E01 | 4.55E00 | 4.3E00 |
| | Std | 4.38E03 | 4.61E-03 | 9.21E-01 | 8.06E00 | 5.8E-01 | 1.84E-01 | 1.38E00 | 2.18E02 | 1.18E-01 | 1.83E-02 |
| | Best | 4.29E00 | 4.31E00 | 4.29E00 | 4.29E00 | 4.52E00 | 4.29E00 | 4.29E00 | 8.24E00 | 4.3E00 | 4.29E00 |
| EP4 | Ave | 2.71E06 | 1.98E00 | 2.11E00 | 1.31E01 | 2.15E00 | 2.03E00 | 2.21E00 | 5.44E01 | 1.98E00 | 1.98E00 |
| | Std | 2.38E06 | 1.03E-15 | 3.11E-01 | 3.42E01 | 1.42E-01 | 3.81E-02 | 3.18E-01 | 5.21E01 | 2.28E-03 | 8.67E-07 |
| | Best | 1.52E02 | 1.98E00 | 1.98E00 | 1.98E00 | 2.04E00 | 1.98E00 | 2E00 | 2.46E00 | 1.98E00 | 1.98E00 |
| EP5 | Ave | 3.08E19 | 1.93E01 | 2.2E01 | 2.43E01 | 3.89E13 | 2.88E01 | 2.72E01 | 1.5E19 | 8.83E14 | 1.85E01 |
| | Std | 2.76E19 | 8.86E-01 | 2.29E00 | 2.99E00 | 1.96E14 | 2.77E00 | 4.27E00 | 7.21E18 | 7.24E14 | 4.76E-02 |
| | Best | 2.42E01 | 1.84E01 | 1.87E01 | 1.88E01 | 2.26E01 | 2.23E01 | 2E01 | 3.12E18 | 1.99E13 | 1.84E01 |
| EP6 | Ave | 7.72E37 | 1.95E04 | 2.01E04 | 3.26E04 | 6.66E24 | 5.8E24 | 3.57E26 | 3.17E26 | 2.37E04 | 1.95E04 |
| | Std | 4.07E38 | 3.78E-09 | 5.22E02 | 1.18E04 | 3.65E25 | 3.18E25 | 1.87E27 | 1.73E27 | 1.72E03 | 1.43E-05 |
| | Best | 1.59E30 | 1.95E04 | 1.95E04 | 1.99E04 | 1.95E04 | 1.96E04 | 2.79E04 | 3.54E04 | 2.08E04 | 1.95E04 |
| EP7 | Ave | 2.67E01 | 2.28E01 | 2.3E01 | 2.31E01 | 2.52E01 | 2.31E01 | 2.31E01 | 2.67E01 | 2.28E01 | 2.28E01 |
| | Std | 1.61E00 | 6.26E-15 | 1.92E-01 | 3.15E-01 | 1.57E00 | 2.96E-01 | 2.36E-01 | 6.24E-01 | 2.18E-03 | 3.81E-09 |
| | Best | 2.44E01 | 2.28E01 | 2.28E01 | 2.28E01 | 2.3E01 | 2.28E01 | 2.28E01 | 2.55E01 | 2.28E01 | 2.28E01 |
| EP8 | Ave | 1.03E01 | 1.34E00 | 1.34E00 | 1.34E00 | 1.36E00 | 1.34E00 | 1.34E00 | 1.38E00 | 1.34E00 | 1.34E00 |
| | Std | 4.4E00 | 4.52E-16 | 1.05E-05 | 2.57E-06 | 1.24E-02 | 3.24E-06 | 3.63E-04 | 1.67E-02 | 4.65E-06 | 1.41E-09 |
| | Best | 1.34E00 | 1.34E00 | 1.34E00 | 1.34E00 | 1.34E00 | 1.34E00 | 1.34E00 | 1.35E00 | 1.34E00 | 1.34E00 |
| EP9 | Ave | 2.7E09 | -8.19E04 | -8.07E04 | -7.73E04 | -8.14E04 | -8.19E04 | -8.15E04 | -7.61E04 | -8.19E04 | -8.19E04 |
| | Std | 9.33E09 | 7.4E-11 | 8.27E02 | 1.01E04 | 1.6E03 | 4.69E-01 | 5.77E02 | 2.84E03 | 7.02E-08 | 7.4E-11 |
| | Best | -8.08E04 | -8.19E04 | -8.19E04 | -8.19E04 | -8.19E04 | -8.19E04 | -8.19E04 | -7.98E04 | -8.19E04 | -8.19E04 |
| EP10 | Ave | 4.16E-01 | 3.14E-01 | 3.14E-01 | 3.18E-01 | 3.64E-01 | 3.16E-01 | 3.14E-01 | 3.23E-01 | 3.14E-01 | 3.14E-01 |
| | Std | 1.25E-02 | 1.13E-16 | 1.13E-16 | 1.58E-02 | 3.86E-02 | 1.43E-02 | 1.13E-16 | 8.58E-03 | 1.13E-16 | 1.13E-16 |
| | Best | 3.14E-01 | 3.14E-01 | 3.14E-01 | 3.14E-01 | 3.14E-01 | 3.14E-01 | 3.14E-01 | 3.14E-01 | 3.14E-01 | 3.14E-01 |
| EP11 | Ave | 3.17E03 | 2.99E03 | 2.99E03 | 3E03 | 3E03 | 3E03 | 3E03 | 3.44E03 | 2.99E03 | 2.99E03 |
| | Std | 1.53E02 | 1.85E-12 | 1.65E-11 | 1.2E01 | 1.75E00 | 2.66E00 | 7.47E00 | 1.58E02 | 1.84E-12 | 1.85E-12 |
| | Best | 3.03E03 | 2.99E03 | 2.99E03 | 2.99E03 | 2.99E03 | 2.99E03 | 2.99E03 | 3.07E03 | 2.99E03 | 2.99E03 |
| EP12 | Ave | 5.97E13 | 1.27E-02 | 1.27E-02 | 1.29E-02 | 1.41E-02 | 1.27E-02 | 1.34E-02 | 1.77E-02 | 1.27E-02 | 1.27E-02 |
| | Std | 8.83E13 | 3.29E-07 | 1.86E-05 | 2.7E-04 | 1.88E-03 | 9.84E-05 | 1.54E-03 | 4.24E-03 | 7.98E-05 | 3.28E-09 |
| | Best | 1.32E-02 | 1.27E-02 | 1.27E-02 | 1.27E-02 | 1.27E-02 | 1.27E-02 | 1.27E-02 | 1.34E-02 | 1.27E-02 | 1.27E-02 |
| EP13 | Ave | 1.35E16 | 1.72E00 | 1.72E00 | 1.8E00 | 3.82E00 | 1.73E00 | 1.9E00 | 2.51E00 | 1.72E00 | 1.72E00 |
| | Std | 5.84E16 | 0E00 | 3.18E-13 | 1.78E-01 | 1.12E00 | 1.58E-03 | 2.46E-01 | 2.04E-01 | 8.04E-06 | 1.11E-09 |
| | Best | 1.81E00 | 1.72E00 | 1.72E00 | 1.72E00 | 2.16E00 | 1.72E00 | 1.72E00 | 2.02E00 | 1.72E00 | 1.72E00 |
| EP14 | Ave | 1.13E06 | 6.18E03 | 6.79E03 | 6.83E03 | 6.88E03 | 6.62E03 | 6.67E03 | 1.16E04 | 6.15E03 | 6.08E03 |
| | Std | 2.48E06 | 1.77E02 | 6.43E02 | 5.67E02 | 4.75E02 | 5.25E02 | 6.44E02 | 2.13E03 | 1.21E02 | 5.67E01 |
| | Best | 6.06E03 | 6.06E03 | 6.06E03 | 6.06E03 | 6.14E03 | 6.06E03 | 6.06E03 | 8.01E03 | 6.06E03 | 6.06E03 |
| Num Of Tops | | 0 | 11 | 6 | 2 | 0 | 4 | 3 | 1 | 9 | **13** |
| Rank | | 9 | 2 | 4 | 7 | 9 | 5 | 6 | 8 | 3 | **1** |

## 3.2 Sensitivity analysis of parameters

There are four parameters which can be adjusted by users when applying AVLA to a specified optimization problem. Among them, the size of the population $N$ and the maximal number of allowed iterations $maxNumIter$ of AVLA are commonly used by nearly all the population-based heuristic algorithms. The size of historical memory $H$ and the maximal number of allowed non-improvement successive iterations $n_R$ are required by AVLA and a few of other heuristics. To investigate the influence of these parameters on the performance of AVLA, we need to choose a special optimization problem which AVLA has not provided with the theoretical optima each time. The problem CF6 in CEC'15 test suite is picked to serve this goal.

In Table 9, we can see SHADE obtains the best resulted mean 6.05E02 for this problem. With the given parameters in Table 2, AVLA gives the second best mean 6.22E02 for this problem. The assumed parameters have the following values: $N$=50, $maxNumIter$=2000, $H$=50, $n_R$=6. Now we will change one parameter's value and fix the others as given to see what will happen to the means resulted from 30 times of independent runs of AVLA to solve CF6 in CEC 2015.

The influence of $N$ is just as we expected the larger the size of population in a reasonable range, the better the result will be. For example, if we raise $N$ to 100, the mean will lower to 6.02E02 with the standard deviation 2.71E00. Similarly, raising the number of allowed maximal iterations will lead to improved result. For example, if we set $maxNumIter$=5000, AVLA will get a lowered mean 6.01E02 with associated standard deviation 2.32E00. Note that the above newly generated results by AVLA are better than that given by SHADE. This just reminds us that we should only accept the comparison conclusion between any algorithms under the preset common condition. If the preset common condition is changed, the comparison conclusion may change accordingly.

To check the influence of varied $H$, we choose a series of values for $H$ and summarize the results associated with these $H$s in Table 12. AVLA obtains the best mean at H=10 and the worst mean at H=100. There is no strict increase or decrease reflected in the data. Based on the above observation, we can see that to choose a proper value of $H$ for a given problem generally needs some trial-and-error. In general, $H$ should not be set bigger than $2H$ because a large $H$ will make the mechanism of success history-based parameter adaption ineffective.

**Table 12** The performance of AVLA with varied $H$.

| $H$ | 5 | 10 | 20 | 30 | 40 | 50 | 60 | 70 | 80 | 90 | 100 |
|---|---|---|---|---|---|---|---|---|---|---|---|
| Ave | 6.09E02 | 6.04E02 | 6.12E02 | 6.11E02 | 6.1E02 | 6.13E02 | 6.16E02 | 6.11E02 | 6.13E02 | 6.15E02 | 6.25E02 |
| Std | 2.19E01 | 7.48E00 | 3E01 | 2.51E01 | 2.29E01 | 3.29E01 | 3.59E01 | 3.14E01 | 3.11E01 | 3.81E01 | 4.69E01 |

With a different $n_R$, AVLA usually obtains a different mean. In Table 13, the results associated with different $n_R$ are listed. From the data given in Table 13, we can see AVLA can obtain a relatively small mean at a low or a high $n_R$. If the value of $n_R$ is relatively small, the reflection by the whole population will be conducted more often; if the value of $n_R$ is relatively big, the reflection by the whole population will be conducted less often. Since when the reflection by the whole population is not conducted at the end of an iteration, the reflection by members with unsatisfied performances will be conducted, $n_R$ can be viewed as a control parameter to balance these two types of reflections. The data in Table 13 shows that a proper value for $n_R$ should be chosen by trial-and error. Most of time the best value of this parameter is problem-specified.

**Table 13** The performance of AVLA with varied $n_R$.

| $n_R$ | 6 | 10 | 20 | 30 | 40 | 50 | 60 | 90 | 100 | 120 | 150 |
|---|---|---|---|---|---|---|---|---|---|---|---|
| Ave | 6.13E02 | 6.21E02 | 6.04E02 | 6.18E02 | 6.22E02 | 6.12E02 | 6.12E02 | 6.14E02 | 6.03E02 | 6.14E02 | 6.09E02 |
| Std | 3.29E01 | 4.43E01 | 4.6E00 | 2.51E01 | 4.56E01 | 3.11E01 | 3.01E01 | 3.11E01 | 3.8E00 | 3.62E01 | 2.38E01 |

## 4. Conclusions

In this paper, we have proposed a novel heuristic algorithm inspired by the learning behaviors of individuals in a group. To verify the efficiency of the new algorithm-AVLA, we apply it to 100 test problems including 29 well-known benchmark problems, CEC 2014 test suite, CEC 2015 learning based test suite, CEC 2022 test suite and 14 engineering optimization problems. By comparing the results obtained by 10 algorithms, we find out that AVLA performs as well as SHADE and both of them outperforms the others in most cases. AVLA and SHADE can be used as a pair of complementary algorithms because they can provide best results for different optimization problems. The variant of AVLA without parameter adaptation called VLA is the best among the left eight algorithms. VLA has its own value due to the fact that it can outperform both AVLA and SHADE in some cases.

There are serval research directions which can be further pursued in the future based on the current understanding of AVLA. By comparing the performances of VLA and AVLA, we know the parameter adaptation is not always helpful when dealing with some types of problems. When and how the parameter adaptation produces effect is an important issue which needs further study. By changing values of the size of memory and the maximal number of non-improvement iterations, we can obtain better results by using AVLA. But the above benefit resulted from adjustment seems problem-specified. In the future, how to add these two parameters in the list of parameters which can be adapted during the optimization should be considered to further improve the performance and resilience of AVLA. Though we have tested AVLA on 100 well-known optimization problems, to implement it to other real-life problems should be conducted continuously in the future to verify its efficiency in a wider area and to find out its weakness for further refinement.

**Disclosure statement**
No potential conflict of interest was reported by the authors.

**Data availability statement**

The data that support the findings of this study are available from the corresponding author upon request.

# Appendix (29 Standard Benchmark Optimization Problems):

**Table 1** Unimodal benchmark functions.

| Functions | Dim | Range | $f_{min}$ |
|---|---|---|---|
| $f_1(x) = \sum_{i=1}^{n} x_i^2$ | 30 | $[-100,100]^n$ | $f(0,0,\ldots,0) = 0$ |
| $f_2(x) = \sum_{i=1}^{n} |x_i| + \prod_{i=1}^{n} |x_i|$ | 30 | $[-10,10]^n$ | $f(0,0,\ldots,0) = 0$ |
| $f_3(x) = \sum_{i=1}^{n} (\sum_{j=1}^{i} x_j)^2$ | 30 | $[-100,100]^n$ | $f(0,0,\ldots,0) = 0$ |
| $f_4(x) = \max_i \{|x_i|, 1 \leq i \leq n\}$ | 30 | $[-100,100]^n$ | $f(0,0,\ldots,0) = 0$ |
| $f_5(x) = \sum_{i=1}^{n-1} [100(x_{i+1} - x_i^2)^2 + (x_i - 1)^2]$ | 30 | $[-30,30]^n$ | $f(1,1,\ldots,1) = 0$ |
| $f_6(x) = \sum_{i=1}^{n} (\lfloor x_i + 0.5 \rfloor)^2$ | 30 | $[-100,100]^n$ | 0 |
| $f_7(x) = \sum_{i=1}^{n} i x_i^4 + random[0,1)$ | 30 | $[-1.28,1.28]^n$ | 0 |

**Table 2** Multimodal benchmark functions.

| Functions | Dim | Range | $f_{min}$ |
|---|---|---|---|
| $f_8(x) = \sum_{i=1}^{n} [-x_i \sin(\sqrt{|x_i|})]$ | 30 | $[-500,500]^n$ | $-418.9829n$ |
| $f_9(x) = \sum_{i=1}^{n} [x_i^2 - 10\cos(2\pi x_i) + 10]$ | 30 | $[-5.12,5.12]^n$ | 0 |
| $f_{10}(x) = -20\exp\left(-0.2\sqrt{\frac{1}{n}\sum_{i=1}^{n} x_i^2}\right) - \exp\left(\frac{1}{n}\sum_{i=1}^{n} \cos(2\pi x_i)\right) + 20 + e$ | 30 | $[-32,32]^n$ | 0 |
| $f_{11}(x) = \sum_{i=1}^{d} \frac{x_i^2}{4000} - \prod_{i=1}^{d} \cos\left(\frac{x_i}{\sqrt{i}}\right) + 1$ | 30 | $[-512,512]^n$ | 0 |
| $f_{12}(x) = \frac{\pi}{n}\{10\sin^2(\pi y_1) + \sum_{i=1}^{n-1}(y_i - 1)^2[1 + 10\sin^2(\pi y_{i+1})] + (y_n - 1)^2\} + \sum_{i=1}^{n} u(x_i, 10, 100, 4)$ | 30 | $[-50,50]^n$ | 0 |
| $y_i = 1 + \frac{x_i+1}{4}, \; u(x_i, a, k, m) = \begin{cases} k(x_i - a)^m & x_i > a \\ 0 & -a \leq x_i \leq a \\ k(-x_i - a)^m & x_i < -a \end{cases}$ | | | |
| $f_{13}(x) = 0.1\{\sin^2(3\pi x_1) + \sum_{i=1}^{n}(x_i - 1)^2[1 + \sin^2(3\pi x_i + 1)] + (x_n - 1)^2[1 + \sin^2(2\pi x_n)]\} + \sum_{i=1}^{n} u(x_i, 5, 100, 4)$ | 30 | $[-50,50]^n$ | 0 |

**Table 3** Fixed dimension multimodal benchmark functions.

| Functions | Dim | Range | $f_{min}$ |
|---|---|---|---|
| $f_{14}(x) = \left(\frac{1}{500} + \sum_{j=1}^{25} \frac{1}{j + \sum_{i=1}^{2}(x_i - a_{ij})^6}\right)^{-1}$ | 2 | $[-65.536, 65.536]^n$ | 1 |
| $f_{15}(x) = \sum_{i=1}^{11} \left[a_i - \frac{x_1(b_i^2 + b_i x_2)}{b_i^2 + b_i x_3 + x_4}\right]^2$ | 4 | $[-5,5]^n$ | 0.0003075 |
| $f_{16}(x) = 4x_1^2 - 2.1x_1^4 + x_1^6/3 + x_1 x_2 - 4x_2^2 + 4x_2^4$ | 2 | $[-5,5]^n$ | -1.0316285 |
| $f_{17}(x) = \left(x_2 - \frac{5.1}{4\pi^2}x_1^2 + \frac{5}{\pi}x_1 - 6\right)^2 + 10\left(1 - \frac{1}{8\pi}\right)\cos x_1 + 10$ | 2 | $[-5,5]^n$ | 0.398 |
| $f_{18}(x) = [1 + (x_1 + x_2 + 1)^2(19 - 14x_1 + 3x_1^2 - 14x_1 + 6x_1 x_2 + 3x_2^2)] \times [30 + (2x_1 - 3x_2)^2(18 - 32x_1 + 12x_1^2 + 48x_2 - 36x_1 x_2 + 27x_2^2)]$ | 2 | $[-2,2]^n$ | 3 |
| $f_{19}(x) = -\sum_{i=1}^{4} c_i \exp\left(-\sum_{j=1}^{3} a_{ij}(x_j - p_{ij})^2\right)$ | 3 | $[0,1]^n$ | -3.86 |
| $f_{20}(x) = -\sum_{i=1}^{4} c_i \exp\left(-\sum_{j=1}^{6} a_{ij}(x_j - p_{ij})^2\right)$ | 6 | $[0,1]^n$ | -3.32 |
| $f_{21}(x) = -\sum_{i=1}^{5} [(x - a_i)(x - a_i)^T + c_i]^{-1}$ | 4 | $[0,10]^n$ | -10.1532 |
| $f_{22}(x) = -\sum_{i=1}^{7} [(x - a_i)(x - a_i)^T + c_i]^{-1}$ | 4 | $[0,10]^n$ | -10.4028 |
| $f_{23}(x) = -\sum_{i=1}^{10} [(x - a_i)(x - a_i)^T + c_i]^{-1}$ | 4 | $[0,10]^n$ | -10.5363 |

**Table 4** Composite benchmark functions.

| Functions | Dim | Range | $f_{min}$ |
|---|---|---|---|
| $F_{24}$ (CF1) <br> $f_1, f_2, f_3, \ldots, f_{10}$ =Sphere Function <br> $[\sigma_1, \sigma_2, \sigma_3, \ldots, \sigma_{10}] = [1,1,1,\ldots,1]$ <br> $[\lambda_1, \lambda_2, \lambda_3, \ldots, \lambda_{10}] = [5/100, 5/100, 5/100, \ldots, 5/100]$ | 10 | [-5,5] | 0 |
| $F_{25}$ (CF2) <br> $f_1, f_2, f_3, \ldots, f_{10}$ =Griewank's Function <br> $[\sigma_1, \sigma_2, \sigma_3, \ldots, \sigma_{10}] = [1,1,1,\ldots,1]$ <br> $[\lambda_1, \lambda_2, \lambda_3, \ldots, \lambda_{10}] = [5/100, 5/100, 5/100, \ldots, 5/100]$ | 10 | [-5,5] | 0 |
| $F_{26}$ (CF3) <br> $f_1, f_2, f_3, \ldots, f_{10}$ =Griewank's Function <br> $[\sigma_1, \sigma_2, \sigma_3, \ldots, \sigma_{10}] = [1,1,1,\ldots,1]$ <br> $[\lambda_1, \lambda_2, \lambda_3, \ldots, \lambda_{10}] = [1,1,1,\ldots,1]$ | 10 | [-5,5] | 0 |
| $F_{27}$ (CF4) <br> $f_1, f_2$= Ackley's Function, $f_3, f_4$= Rastrigin's Function, <br> $f_5, f_6$= Weierstrass Function, $f_7, f_8$= Griewank's Function, <br> $f_9, f_{10}$= Sphere Function <br> $[\sigma_1, \sigma_2, \sigma_3, \ldots, \sigma_{10}] = [1,1,1,\ldots,1]$ <br> $[\lambda_1, \lambda_2, \lambda_3, \ldots, \lambda_{10}]$ =[5/32, 5/32, 1, 1, 5/0.5, 5/0.5, 5/100, 5/100, 5/100, 5/100] | 10 | [-5,5] | 0 |
| $F_{28}$ (CF5) <br> $f_1, f_2$= Rastrigin s Function, $f_3, f_4$= Weierstrass Function, <br> $f_5, f_6$= Griewank s Function, $f_7, f_8$= Ackley s Function, <br> $f_9, f_{10}$ = Sphere Function <br> $[\sigma_1, \sigma_2, \sigma_3, \ldots, \sigma_{10}] = [1,1,1,\ldots,1]$ <br> $[\lambda_1, \lambda_2, \lambda_3, \ldots, \lambda_{10}]$ =[1/5, 1/5, 5/0.5, 5/0.5, 5/100, 5/100, 5/32, 5/32, 5/100, 5/100] | 10 | [-5,5] | 0 |
| $f_{29}$ (CF6) <br> $f_1, f_2$= Rastrigin's Function, $f_3, f_4$= Weierstrass' Function, <br> $f_5, f_6$= Griewank's Function, $f_7, f_8$= Ackley's Function, <br> $f_9, f_{10}$= Sphere Function <br> $[\sigma_1, \sigma_2, \sigma_3, \ldots, \sigma_{10}]$ =[0.1, 0.2, 0.3, 0.4, 0.5, 0.6, 0.7, 0.8, 0.9, 1] <br> $[\lambda_1, \lambda_2, \lambda_3, \ldots, \lambda_{10}]$ =[0.1 ∗ 1/5, 0.2 ∗ 1/5, 0.3 ∗ 5/0.5, 0.4 ∗ 5/0.5, 0.5 ∗ 5/100, 0.6 ∗ 5/100, 0.7 ∗ 5/32, 0.8 ∗ 5/32, 0.9 ∗ 5/100, 1 ∗ 5/100] | 10 | [-5,5] | 0 |